\documentclass{article}
\usepackage{blindtext}
\usepackage{geometry}
\usepackage{amssymb}
\usepackage{makecell}

\usepackage{subfiles} 
\usepackage{graphicx} 
\usepackage[utf8]{inputenc}
\usepackage[english]{babel}  
\usepackage[T1]{fontenc}
\usepackage{mathtools}
\usepackage{amsthm}
\usepackage{amssymb}
\usepackage{amsmath}
\usepackage{fourier}
\usepackage[dvipsnames]{xcolor}
\usepackage{comment}
\usepackage{authblk}
\usepackage[colorlinks=true,linkcolor=black,anchorcolor=black,citecolor=black,filecolor=black,menucolor=black,runcolor=black,urlcolor=black]{hyperref} 

\theoremstyle{definition}
\newtheorem{definition}{Definition}[section]
\newtheorem{theorem}{Lemma}[section]

\newcommand{\R}{\mathbb{R}}
\newcommand{\Z}{\mathbb{Z}}
\newcommand{\N}{\mathbb{N}}
\newcommand{\C}{\mathbb{C}}
\newcommand{\quat}{\mathbb{H}}
\newcommand{\W}{\mathcal{W}}
\newcommand{\hil}{\mathcal{H}}
\newcommand{\Lim}[1]{\raisebox{0.5ex}{\scalebox{0.8}{$\displaystyle \lim_{#1}\;$}}}

\newtheorem{proposition}{Proposition}[section]

\title{Comparison between two approaches to classify topological insulators using $K$-theory}
\date{January 2024}
\author{Lorenzo Scaglione}
\affil{Institut Camille Jordan, Université Claude Bernard Lyon 1, 43 boulevard du 11 novembre 1918, F-69622 Villeurbanne Cedex}

\begin{document}

\maketitle
\begin{abstract}
    We compare two approaches which use $K$-theory for $C^*$-algebras to classify symmetry protected topological phases of quantum systems described in the one particle approximation. In the approach by Kellendonk, which is more abstract and more general, the algebra remains unspecified and the symmetries are defined using gradings and real structures. In the approach by Alldridge et al., the algebra is physically motivated and the symmetries implemented by generators which commute with the Hamiltonian. Both approaches use van Daele’s version of $K$-theory. We show that the second approach is a special case of the first one. We highlight the role played by two of the symmetries: charge conservation and spin rotation symmetry.
\\
    
Acknowledgement: We thank Christopher Max for sharing his insight into this problem and my supervisor Johannes Kellendonk for fruitful discussions.

\end{abstract}
\tableofcontents
\newpage

\section{Introduction}
Topological insulators are materials with very interesting properties from a physical point of view: they are electrical insulators in the bulk, while on the surface they can be characterised by the presence of robust currents. These properties are specific to the \textit{topological phase} of the material, i.e. their equivalence class with respect to continuous deformations in a given topology. This means that small external perturbations or microscopic irregularities of the material, which we can consider as continuous deformations of the insulator, will not affect the intrinsic physical properties of the phase. To treat this problem mathematically, first of all we identify an insulator with the Hamiltonian which describes the behaviour of electrons in the material. As usual in physics, Hamiltonians of insulators are self-adjoint operators over an Hilbert space $\hil$ which have a gap in their spectrum at the Fermi energy, which we ﬁx to be at $0$ (this is equivalent to say that it that the operator is invertible). However, equipped with the norm topology $\mathcal{B}(\hil)$ is non-separable (if the dimension is inﬁnite) and its $K$-theory trivial. For this reason we choose as the relevant topological space a sub-$C^*$-algebra of $\mathcal{B}(\hil)$ equipped the norm topology.
We could impose extra physical symmetries to the system and our Hamiltonians should satisfy some extra constraints.

The mathematical problem we want to solve is to give a classification of the topological phases, namely to describe path-connected components of the set of Hamiltonians associated with the insulators. More precisely, we consider a $C^*$-algebra $A$ (the choice of this algebra comes from the specific physical situation) and a set of symmetries $\mathcal{R}$. The general goal is to describe the set of Hamiltonians satisfying the symmetries of $\mathcal{R}$ up to homotopy (therefore the set of topological phases) denoted by
$$GL^{s.a.}_{h}(A, \mathcal{R}).$$

Here we consider two approaches, \cite{kellendonk} and \cite{zirnbauer}, which use $K$-theory for $C^*$-algebras to classify topological phases when symmetries are imposed. The first is more abstract: in \cite{kellendonk}, the algebra containing the abstract insulators is arbitrary and the symmetries are defined abstractly using automorphisms and real structures. We will see that it is useful to have a reference real structures in order to bring back the study of all other possible real structures to a common real structure. In \cite{zirnbauer}, starting from a \textit{non-interacting free fermion} description of the solid, we obtain a  very precise and physically well motivated algebra and symmetries. In this context, a reference real structure naturally appears from the beginning of the construction of the algebra: therefore, we deal with a Real $C^*$-algebra (a complex $C^*$-algebra endowed with a real structure) and this structure will be the natural candidate to study all additional symmetries. 

The comparison of the approaches is possible because both approaches use a common tool to describe the set $GL^{s.a.}_{h}(A, \mathcal{R})$, namely Van Daele's version of $K$-theory. We will show Van Daele's construction and we will see why it is surprisingly natural to describe topological phases. Applying this construction to the relevant $C^*$-algebra for a choice of symmetries, we obtain an Abelian group. We will often write the symbol 
$$GL^{s.a.}_{h}(A, \mathcal{R})\leadsto DK(A', \alpha, \mathfrak{r})$$
to say that we associate to $GL^{s.a.}_{h}(A, \mathcal{R})$ the Van Daele's Abelian group $DK(A', \alpha, \mathfrak{r})$. These two objects are not exactly the same, but we will see that $DK(A', \alpha, \mathfrak{r})$ can give important information about $GL^{s.a.}_{h}(A, \mathcal{R})$ and that, physically, this association is well-justified. An important remark is however that the classification given by Van Daele's construction is relative: indeed, to obtain a Van Daele's group we have to choose a reference element. Two different choices lead to isomorphic Van Daele's groups, but the isomorphism is not canonical.   

Our contribution to the general goal presented above is to establish a link between the results obtained with the two approaches and, in particular, to show how the concrete approach is a special case of a more general one. We would like all of this to be coherent with the \textit{Tenfold Way} of the Kitaev table \cite{kitaev2009periodic} in order to compare these approach to the classical literature.  
More precisely, in the \cite{zirnbauer}, approach we have ten cases to deal with, which correspond to the possible combinations of symmetries: we have two complex cases and eight real cases. For each case, we can show how to express the symmetries in the abstract language and check that in both cases we obtain the same $K$-group. What is interesting in doing so is that we will understand the essence of symmetries and which is their mathematical role. The main merit of this comparison is to emphasise the role of two type of symmetries: the charge conservation and the Spin Rotations Symmetry. Charge conservation has the role of reducing the algebra. Indeed, in this case we can just consider a smaller complex sub-algebra without loss of information. Spin Rotations Symmetry allows us to factor the algebra using quaternions; therefore, we can consider the additional symmetries (if any) over a reduced algebra. We will see that considering the quaternionic part will change the relative sign of the symmetries (an even Time Reversal Symmetry becomes an odd Time Reversal Symmetry, for example) and implying a shift of four units of the corresponding $K$-group index.    

This work is organised as follows. In the first part, we will describe the abstract approach and in particular Van Daele's version of $K$-theory. In the second part, we will deal with the physical construction of a possible $C^*$-algebra describing topological insulators and show the table of $K$-groups obtained in \cite{zirnbauer}. Finally, the core of our work is the third part: we will establish an explicit relation between the two approaches and show case by case that we obtain the same $K$-groups. In the final paragraph we will show the table summarising the translation from the abstract approach to the concrete one. 

\section{Kellendonk's approach}
Here, an \textit{abstract} insulator is a self-adjoint and invertible element of an arbitrary $C^*$-algebra. Symmetries are described by a linear automorphism (in this case, we have a chiral symmetry) and/or an antilinear automorphism (we have a Time Reversal Symmetry or Particle Hole Symmetry) of order $2$. Compared with the concrete approach, here we have much more freedom over the choice of automorphisms, but we can decide to impose constraints to refine the results: for example, we can make the hypothesis of a particular form for the automorphisms (we will often meet \textit{inner} chiral symmetries) or of a particular relationship between symmetries. It turns out that these additional constraints are satisfied in the concrete case which therefore fits the abstract framework .

\subsection{Symmetries}
\label{sec:physical_symmetries}
In the following, we will consider Hamiltonians which could possibly satisfy one or more symmetry conditions. The physical symmetries are translated into algebraic conditions with the help of order two linear or anti-linear $*$-automorphism. We give the following definitions.

\begin{definition}
	A \textit{grading} on a (complex or real) $C^*$-algebra $A$ is a
	$*$-automorphism $\gamma$ of order two ($\gamma ^2=\text{id}$).
\end{definition}
\begin{definition}
	A \textit{real structure} $\gamma$ on a complex vector space is an anti-linear ($\gamma(\text{i}\cdot)=-\text{i}\gamma(\cdot)$) automorphism of order two (i.e. $\gamma^2=\text{id}$).
\end{definition}
\begin{definition}
	A \textit{real structure} on a complex graded $C^*$-algebra $(A, \gamma)$ is a real structure $\mathfrak{r}$ on the vector space $A$ which is a $*$-automorphism of $A$ as $C^*$-algebra, commutes with the grading $\gamma$ and preserves the norm.
\end{definition}
\begin{definition}\label{RealC*alg}
	A complex $C^*$-algebra with a real structure $(A, \mathfrak{r})$ is called a \textit{Real $C^*$-algebra} or a $C^{*,r}$-algebra. The $\mathfrak{r}$-invariant elements furnish a real $C^*$-algebra which we call the real sub-algebra of the $C^{*,r}$-algebra and note it with $A^{\mathfrak{r}}$.
\end{definition}
A graded $C^*$-algebra $(A, \gamma)$ is a $C^*$-algebra $A$ equipped with a grading $\gamma$. $\gamma$-invariant elements are called \textit{even}, and elements with $\gamma(a) = -a$ are called \textit{odd}.

In the following we will often consider special gradings, inner gradings and balanced gradings.
\begin{definition}
A grading $\gamma$ is called \textit{inner} if $\gamma = \text{Ad}_\Gamma$ for some self-adjoint unitary $\Gamma$ in $\mathcal{M}(A)$, the multiplier algebra of $A$. The self-adjoint unitary $\Gamma$ is called the generator of $\gamma$ or the grading operator.
\end{definition}

\begin{definition}
A grading on a $C^*$-algebra is called \textit{balanced} if the $C^*$-algebra contains an odd self-adjoint unitary $e$.
\end{definition}
In our physical setting, the hypothesis of having a balanced grading is equivalent to say that it exists at least one insulator.

We have three relevant symmetries.
\begin{definition}
	Let $A$ be a $C^*$-algebra and $h \in A$ an abstract Hamiltonian.
	\begin{itemize}
		\item $h$ has \textit{chiral} symmetry if there is a grading $\gamma$ on $A$ such that $\gamma(h) = -h$;
		\item $h$ has \textit{Time Reversal Symmetry (TRS)} if there is a real structure $\mathfrak{t}$ on $A$ such that $\mathfrak{t}(h) = h$;
		\item $h$ has \textit{Particle Hole Symmetry (PHS)} if there is a real structure $\mathfrak{p}$ on $A$ such that $\mathfrak{p}(h) = -h$.
	\end{itemize}
\end{definition}
We have that the product of a Time Reversal Symmetry together with a Particle Hole Symmetry yields a chiral symmetry. We thus
have the following combinations:
\begin{itemize}
	\item no symmetry
	\item chiral symmetry
	\item Time Reversal Symmetry
	\item Particle Hole Symmetry
	\item chiral symmetry and Time Reversal Symmetry (and therefore, also Particle Hole Symmetry).     
\end{itemize}
These five cases will be refined (we will define an odd/even TRS and PHS) and we will obtain two complex cases and eight real cases.

\subsection{Clifford algebras}
\label{sec:Clifford}
Clifford algebras are a technical tool which will allow us to give us a unified framework to treat all cases of symmetries (the five cases of symmetries presented at the end of section \ref{sec:physical_symmetries}).
\par 

We give a definition of Clifford algebra suited for our case.
\begin{definition}
The Clifford algebra $Cl_{r,s}$ is the graded real $C^*$-algebra generated by $r$ self-adjoint generators $e_1,\dots, e_r$ which square to $1$ and $s$ anti-self-adjoint generators $f_1 ,\dots, f_s$ which square to $-1$ and all generators anti-commute pairwise. The grading is deﬁned by declaring the generators to be odd.
\end{definition}
\begin{definition}
The complex Clifford algebra $\C l_{r+s}$ is the complexification of $Cl_{r,s}$, $\C l_{r+s}=\C \otimes Cl_{r,s}$.
\end{definition}
The generators define completely the grading, which is noted st. We remind the expression of the Pauli matrices
$$\sigma_x=\begin{pmatrix}
0 & 1\\
1 & 0
\end{pmatrix}
\quad
\sigma_y=\begin{pmatrix}
0 & -\text{i}\\
\text{i} & 0
\end{pmatrix}
\quad 
\sigma_z=\begin{pmatrix}
1 & 0\\
0 & -1
\end{pmatrix}
$$
Recall that the quaternions $\mathbb{H}$ form a real $C^*$-algebra which is spanned as a real vector space by $\{ 1, \text{i} \sigma_x, \text{i} \sigma_y, \text{i} \sigma_z \}$. We consider $\mathbb{H}$ always as trivially graded.

The graded tensor product between two graded algebras will be written $\hat{\otimes}$. It can be understood as the ordinary tensor product as far as the linear structure is concerned, but multiplication and involution take care of the Koszul sign rule:
$$(a_1\hat{\otimes}b_1)(a_2\hat{\otimes}b_2)=(-1)^{|a_2||b_1|}(a_1 a_2\hat{\otimes}b_1 b_2) \quad (a\hat{\otimes}b)^*=(-1)^{|a||b|}(a^*\hat{\otimes}b^*)$$
We have the following useful isomorphism (see \cite{kellendonk}, paragraph 3.1.2)
\begin{equation}\label{Cl1}
    \left( Cl_{r,s}\hat{\otimes}\C l_{r',s'},\text{st}\otimes\text{st}\right)\cong \left(Cl_{r+r',s+s'},\text{st}\right)
\end{equation}
\begin{equation}\label{Cl2}
    \left( \quat\otimes\C l_{0,1},\text{id}\otimes\text{st}\right)\cong \left(Cl_{0,4},\text{st}\right)
\end{equation}

\subsection{Van Daele's approach to $K$-theory}
\label{sec:VD_K-theory}
In \cite{VD}, A. Van Daele proposes a new approach to $K$-theory for Banach algebras (not just for $C^*$-algebras, but we will apply this theory just to the $C^*$-algebra case): here, we want to associate an abelian group to the homotopy classes of a subset of the elements of the Banach algebra which satisfy a symmetry condition. This is exactly what we want to obtain physically and Van Daele's approach applies directly to our case.
Finally, in some cases, Van Daele manages to link (i.e. find an isomorphism between) his $K$-groups and the standard complex and real $K$-groups.
\par 
We describe here briefly the construction. Let $A$ be a real or complex $C^*$-algebra with a grading $\alpha$.
As a first step, we just consider the case of an unital $C^*$-algebra.
We define
$$
\mathcal{F}(A, \alpha )=\left \{ a \in A \mid a \text{ self-adjoint and unitary, } \alpha (a)=-a \right \}
$$
We first assume that $\mathcal{F}(A, \alpha )$ is not empty, i.e. the grading is balanced (physically it means that we have at least one insulator).
We denote by $F(A, \alpha)$ the set of homotopy equivalence classes of $\mathcal{F}(A, \alpha)$.\\
We let
$$
\mathcal{F}_{n}(A, \alpha)=\mathcal{F}\left(M_{n}(A), \alpha_n\right)\quad \text { and }\quad F_{n}(A)=F\left(M_{n}(A), \alpha_n\right)
$$
where $\alpha_n$ is application of $\alpha$ element-by-element on $M_n(A)$.
Analogously to the construction of $K_0(A)$ and $K_1(A)$, we define the direct sum $(x,y)\mapsto x\oplus y$ as a map from $\mathcal{F}_n(A)\times \mathcal{F}_m(A)$ to $\mathcal{F}_{m+n}(A)$ such that 
$$(x,y)\mapsto x\oplus y=\begin{pmatrix}
x & 0\\
0 & y
\end{pmatrix}.$$
\begin{definition}
Choose an element $e\in \mathcal{F}(A, \alpha )$ and define $DK_e(A)=\Lim{\longrightarrow} F_n(A, \alpha)$ where the inductive limit is taken with respect to the maps $x\rightarrow x\oplus e$ from $\mathcal{F}_n(A, \alpha )$ to $\mathcal{F}_{n+1}(A, \alpha )$. We will denote by $[x]$ the image of $x$ in the inductive limit.
\end{definition}
The definition of a direct sum on $\mathcal{F}_n(A)\times \mathcal{F}_m(A)$ induces a direct sum also on $F_n(A)\times F_m(A)$ and $DK_e(A)\times DK_e(A)$.
It is interesting that we have to choose a reference element $e$: physically, it means that we must choose a reference insulator.

We have the following propositions.
\begin{proposition}
The direct sum induces an abelian semi-group structure on $DK_e(A)$. The neutral element is $[e]$.
\end{proposition}
\begin{proposition}
If $e$ and $-e$ are homotopic in $\mathcal{F}(A, \alpha)$ then $DK_e(A)$ is a group.
\end{proposition}
Now we define the Van Daele $K$-group $DK(A)$ for an unital $C^*$-algebra where $\mathcal{F}(A, \alpha )$ could be possibly empty.
\begin{definition}
    Let $A$ be a unital $C^*$-algebra. Let $DK(A)$ be the group $DK_e(M_4(A))$ where on $M_4(A)$ the automorphism $\gamma$ is given by the element-wise application of the automorphism $\gamma'$ on $M_2(A)$ where
    $$\gamma'\left(\begin{pmatrix}
        a & b\\
        c & d
    \end{pmatrix}\right)=
    \begin{pmatrix}
        \alpha(a) & -\alpha(b)\\
        -\alpha(c) & \alpha(d)
    \end{pmatrix}$$
    and 
    $$e=\begin{pmatrix}
        f & 0\\
        0 & -f
    \end{pmatrix}\quad \text{where}\quad f=\begin{pmatrix}
        0 & 1\\
        1 & 0
    \end{pmatrix}$$
\end{definition}
We can verify that the definition is well-posed. Indeed, we built an element $e\in \mathcal{F}(M_4(A), \alpha )$ (therefore $DK_e(M_4(A))$ is defined as in the previous construction) which is homotopic to $-e$ in $\mathcal{F}(M_4(A), \alpha )$ ($DK_e(M_4(A))$ is a group). The following proposition establishes a link between $DK(A)$ and $DK_{e_0}(A)$ in the case where $\mathcal{F}(A, \alpha )$ is not empty.

\begin{proposition}
    If $A$ has an element $e_0\in \mathcal{F}(A, \alpha)$ then $DK(A)$ is the Grothendieck group of the semi-group $DK_{e_0}(A)$.
\end{proposition}
Finally, we drop the requirement that $A$ has an identity. We let $A^+ = \{(a, \lambda) : a\in A, \lambda\in \C\}$ and we equip $A^+$ with the usual $C^*$-algebra structure. We denote the identity in $A^+$ by $1$ and identify $A$ with its canonical image in $A^+$. So $A^+=A+\C 1$. Let $\phi$ be the canonical homomorphism from $A^+$ to $\C$. We extend the grading $\alpha$ from $A$ to $A^+$ in the obvious unique way saying that the element $1$ is even and we still use $\alpha$ for this extension. It follows that $\phi$ induces a homomorphism $\phi_*:DK(A^+)\rightarrow DK(\C)$.
\begin{definition}
    We define $\Tilde{DK}(A)$ as the kernel of the map $\phi_*:DK(A^+)\rightarrow DK(\C)$. In particular we know that $DK(\C)=\{0\}$ so that $\Tilde{DK}(A) = DK(A^+)$.
\end{definition}
We have the following result.
\begin{proposition}
    If A has an identity then $\Tilde{DK}(A)$ and $DK(A)$ are isomorphic.
\end{proposition}
Therefore, for an arbitrary $C^*$-algebra $A$ we can define its Van Daele's group as $\Tilde{DK}(A)$ (which we still denote $DK(A)$). However, when $\mathcal{F}(A, \alpha )$ is not empty it is better to work with the group $DK_e(A)$ for the sake of interpretation (having in mind that the isomorphism between $DK_e(A)$ and $DK(A)$ depends on the choice of $e$). 

\par
Before introducing Van Daele's higher $K$-groups, we give a useful result. We have
\begin{equation}\label{isomorphism_via_reference_elem}
    D K_{e}\left(A \hat{\otimes} C l_{r, s}, \gamma \otimes \mathrm{st}\right) \cong D K_{e}\left(A \hat{\otimes} C l_{r+1, s+1}, \gamma \otimes \mathrm{st}\right)
\end{equation}
In particular, it implies that $D K_{e}\left(A \hat{\otimes} C l_{r, s}, \gamma \otimes \mathrm{st}\right)$ depends only on the difference $s-r$ due to stabilisation.
We present here Van Daele's higher $K$-groups.
\begin{definition}
Let $(A, \gamma)$ be a unital graded complex or real $C^*$-algebra. The $K_ n$-group of $(A, \gamma)$ in Van Daele’s
formulation is

$$
K_{n}(A, \gamma)\coloneqq D K \left(A \hat{\otimes} C l_{1, n}, \gamma \otimes \mathrm{st}\right).
$$
For $C^{*,r}$-algebras, it is defined as $K_n(A, \gamma, \mathfrak{r})=K_n(A^\mathfrak{r},\gamma)$
\end{definition}
A fundamental property of the $K$-groups $K_{n}(A, \gamma)$ is their periodicity. We have $K_n\cong K_{n-8}$ and if we are in the complex case we also have $K_n\cong K_{n-2}$ (we use here the isomorphism (\ref{isomorphism_via_reference_elem}) and therefore these latter isomorphisms depend on the choice of the reference element).

We conclude this section with a list of useful results in some particular cases.
\begin{itemize}
    \item If the algebra is trivially graded:
    \begin{enumerate}
    
        \item case of a complex $C^*$-algebra $(A,\text{id})$\\
        Van Daele shows that
        \begin{equation}\label{KUi}
            K_{n}(A, \mathrm{id})=D K\left(A \otimes \mathbb{C} l_{n+1}, \mathrm{id} \otimes \mathrm{st}\right) \cong K U_{n}(A)
        \end{equation}
where $KU_n (A)$ is the standard $K_n$-group ($n$ taken modulo $2$) of $A$ seen as ungraded complex $C^*$-algebra.        
        \item case of a $C^{*,r}$-algebra $(A,\text{id}, \mathfrak{r})$\\
        We have
        \begin{equation}\label{KOi}
            D K\left(A \otimes \mathbb{C} l_{r+s}, \mathrm{id} \otimes \mathrm{st}, \mathfrak{r} \otimes \mathfrak{l}_{r, s}\right) \cong K_{s-r+1}\left(A^{\mathfrak{r}}, \mathrm{id}\right) \cong K O_{s-r+1}\left(A^{\mathfrak{r}}\right)
        \end{equation}
where $\mathfrak{l}_{r, s}$ is the real structure on $\C l_{r+s}$ such that $\C l_{r+s}^{\mathfrak{l}_{r, s}}=Cl_{r, s}$ and $KO_n (A)$ is the $K_n$-group of $A$ seen as real ungraded $C^*$-algebra.
    \end{enumerate}
    
    \item If the grading of the algebra is inner:
    \begin{enumerate}
        \item case of a complex $C^*$-algebra $(A,\gamma)$\\
        We have 

        \begin{equation}\label{KU1}
            D K(A, \gamma) \cong K_{1}(A, \mathrm{id}) \cong K U_{1}(A)
        \end{equation}
        \item case of a $C^{*,r}$-algebra $(A,\gamma, \mathfrak{r})$, where $\gamma=\text{Ad}_\Gamma$\\
        If $\mathfrak{r}(\Gamma)=\Gamma$, then
        \begin{equation} \label{TRS_realGrading}
            D K(A, \gamma, \mathfrak{r}) \cong K_{1}\left(A^{\mathfrak{r}}, \mathrm{id}\right) \cong K O_{1}\left(A^{\mathfrak{r}}\right)
        \end{equation}
        
        If $\mathfrak{r}(\Gamma)=-\Gamma$, then
\begin{equation}\label{TRS_immGrading}
    D K(A, \gamma, \mathfrak{r}) \cong K_{-1}\left(A^{\mathfrak{r}}, \mathrm{id}\right) \cong K O_{-1}\left(A^{\mathfrak{r}}\right)
\end{equation}
    \end{enumerate}
    
\end{itemize}

\subsection{How to use Van Daele's construction}
\label{sec:VD}

The goal is to classify symmetric Hamiltonians, namely self-adjoint invertible elements of a $C^*$-algebra which satisfy one or more symmetry conditions. Van Daele's construction (presented in \ref{sec:VD_K-theory}) allows us to obtain a $K$-group corresponding to the homotopy classes in the set of unitary self-adjoint elements which satisfy a chiral symmetry. We will see how to adapt each case in order to apply Van Daele's construction. 
\par 
In \cite{kellendonk}, the author associates to each symmetry case a $K$-group. Two classifications are presented. The first one is a rough classification, under minimal assumptions. The second one is finer and it is a classification with respect to a reference real structure. The classifications can be even more refined if we assume the grading to be inner. 
\par 
Before beginning the description of the classification, we give a general result which is fundamental to apply Van Daele's construction to all cases of symmetry.
In a $C^*$-algebra $A$, any invertible self-adjoint element $h$ is homotopic to a self-adjoint unitary, namely its sign $\text{sgn}(h)=h\lvert h\rvert^{-1}$; this is referred to as \textit{spectral ﬂattening}. 
Equivalently, we can introduce $p_F(h)$, the \textit{Fermi projector} of $h$, and define $\text{sgn}(h)$ as
$$\text{sgn}(h)=p_F(h)-p_F^\perp (h).$$
It can be shown that the $\text{sgn}$ function induces a bijection between the homotopy class of odd self-adjoint invertible elements and the homotopy class of odd unitary self-adjoint elements. In this way, we are allowed to consider just odd unitary self-adjoint elements. 
\par
We will show now how to adapt the Van Daele's construction to each case of symmetry. We consider the graded $(A, \gamma)$, with, for simplicity, $A$ unital. Most of the following results are still valid in the non-unital case.
\subsubsection{No symmetry}\label{noSym}

No symmetry is equivalent to a trivial grading, $\gamma=\text{id}$, as in this case every element is symmetric. We want to describe an insulator without symmetry within the framework of Van Daele. The trick to recover an odd element is to extend the algebra $A$ to $A\otimes \C l_1$.

To see it, first of all we consider the isomorphism between $(\C l_1,\text{st})$ and $(\C\oplus \C,\phi)$, which maps the generator $e_1$ of $\C l_1$ to $(1,-1)$ and where $\phi$ is the "flip" isomorphism. Then, the map $\left(x,-x\right)\mapsto x$ from self-adjoint odd unitaries of $\left( A\otimes \C l_1, \text{id}\otimes \text{st}\right)$ to self-adjoint unitaries of $\left( A,\text{id}\right)$ is a bijection. The advantage is that now each insulator is associated to an odd element of a $C^*$-algebra and we can apply straightforwardly Van Daele's construction. Thanks to the isomorphism \ref{KUi}, we get the $K$-group
$$
D K\left(A \otimes \mathbb{C} l_{1}, \mathrm{id} \otimes \text{st}\right)=K_{0}(A, \mathrm{id}) \cong K U_{0}(A)
$$
\subsubsection{Chiral symmetry}\label{chiralSym}
Now the grading $\gamma$ is non trivial and we assume that the grading is balanced, i.e. it exists at least one self-adjoint odd element (otherwise, the grading is not physically very interesting). We can directly use Van Daele's construction, which gives the $K$-group
$$DK(A,\gamma)=K_1(A, \gamma)$$
If the grading $\gamma$ is inner, then isomorphism (\ref{KU1}) gives
$$DK(A,\gamma)=KU_1(A)$$

\subsubsection{Time Reversal Symmetry}\label{TRSSym}
We consider the trivially graded $C^{*,r}$-algebra $(A, \mathfrak{t})$. Analogously to the case without symmetry, we tensor $A$ with $(\C l_1, \text{st},\mathfrak{l}_{1,0})\cong (\C \oplus \C, \phi,\mathfrak{c})$ (where $\phi$ is the "flip" isomorphism and $\mathfrak{c}$ is the complex conjugation) to recover odd elements; indeed, the map

\begin{align*}
  p\colon (A\otimes (\C \oplus \C), \text{id}\otimes \phi,\mathfrak{t}\otimes \mathfrak{c} ) & \longrightarrow (A, \mathfrak{t}) \\
  x\otimes (1,-1) & \longmapsto x 
  \end{align*}
is a bijection between real odd self-adjoint elements of $A\otimes \C l_1$ and self-adjoint real elements of $A$. We can see it if we remark that real elements of $A\otimes \C l_1$ are $\R$-linear combinations of elements $a\otimes b$ with $a\in A$ and $b\in \C \oplus \C$ real. The other possibility would be $a'\in A$ and $b'\in \C \oplus \C$ imaginary. But then $a'=\text{i}a$ and $b'=\text{i}b$ for some real $a$ and $b$ and therefore $a'\otimes b'=-(a\otimes b)$ is again a $\R$-linear combinations of elements $a\otimes b$ with $a\in A$ and $b\in \C \oplus \C$ real. Therefore, the previous definition characterises $p$ on all real odd elements ($(1,-1)$ is the generator of real odd elements of $\C \oplus \C$).
Therefore, we get the $K$-group (isomorphism \ref{KOi})
\begin{equation*}
    D K\left(A \otimes \mathbb{C} l_{1}, \mathrm{id} \otimes \text{st}, \mathfrak{t} \otimes \mathfrak{l}_{1,0}\right)=K_{0}\left(A^{\mathfrak{t}}, \mathrm{id}\right) \cong K O_{0}\left(A^{\mathfrak{t}}\right)
\end{equation*}

\subsubsection{Particle Hole Symmetry}\label{PHSSym}
We consider the trivially graded $C^{*,r}$-algebra $(A, \mathfrak{p})$. Now we consider the algebra 
$$\left(A \otimes \mathbb{C} l_{1}, \mathrm{id} \otimes \text{st}, \mathfrak{p}\otimes \mathfrak{l}_{0,1}\right)\cong \left(A \otimes (\C\oplus \C), \mathrm{id} \otimes \phi, \mathfrak{p}\otimes\phi\circ\mathfrak{c}\right)$$
Similarly to the previous case, we see that there is a bijection between the real odd elements of $A \otimes \mathbb{C} l_{1}$ and the imaginary elements of $A$ (where the map is $h\otimes \left(1,-1\right)\mapsto h$).
We get the $K$-group
\begin{equation}\label{PHS}
    \begin{aligned}
D K\left(A \otimes \mathbb{C} l_{1}, \mathrm{id} \otimes \text{st}, \mathfrak{p} \otimes \mathfrak{l}_{0,1}\right) &=D K\left(A^{\mathfrak{p}} \otimes C l_{0,1}, \mathrm{id} \otimes \text{st}\right)=K_{2}\left(A^{\mathfrak{p}}, \mathrm{id}\right) \\
& \cong K O_{2}\left(A^{\mathfrak{p}}\right) .
\end{aligned}
\end{equation}
\subsubsection{Chiral symmetry and Time Reversal Symmetry}\label{chiralTRSSym}
We have a real structure $\mathfrak{t}$ and the grading $\gamma$ is balanced.
We get immediately the $K$-group
$$DK(A,\gamma, \mathfrak{t})=K_1(A^\mathfrak{t}, \gamma)$$
\par 
If the grading is inner, we get a finer classification. Isomorphisms (\ref{TRS_realGrading}) and (\ref{TRS_immGrading}) give us

\begin{equation}\label{inner_chiral}
    \begin{array}{cc} 
\hline \text { Symmetries }  & K \text {-group } \\
\hline 
\text { Real inner chiral }  & K O_{1}\left(A^{\mathfrak{t}}\right) \\
\text { Imag. inner chiral }  & K O_{-1}\left(A^{\mathfrak{t}}\right) \\
\hline
\end{array}
\end{equation}

\subsubsection{Classification w.r.t. a reference real structure}
We would like to compare $K O_{i}\left(A^{\mathfrak{t}}\right)$ and $K O_{i}\left(A^{\mathfrak{p}}\right)$. For this reason, we choose a reference real structure $\mathfrak{f}$ over the $C^*$-algebra. 
We define the concept of \textit{inner related} and \textit{inner conjugated} real structures.
\begin{definition}
Let $A$ be a graded complex $C^*$-algebra with two real structures $\mathfrak{r}$, $\mathfrak{s}$. We call them \textit{inner related} if there exists a unitary $u$ (the generator of $\mathfrak{s} \circ \mathfrak{r}$) in the multiplier algebra of $A$ such that $\mathfrak{s} \circ \mathfrak{r} = \text{Ad}_u$.
\end{definition}
\begin{definition}
\label{inner_conj}
Let $A$ be a graded complex $C^*$-algebra with two real structures $\mathfrak{r}$, $\mathfrak{s}$. We call them \textit{inner conjugate} if there exists a unitary $w$ in the multiplier algebra of $A$ such that $\mathfrak{s}\circ \text{Ad}_w = \text{Ad}_w \circ \mathfrak{r}$ and $\text{Ad}_w$ preserves the grading.
\end{definition}
Inner conjugation is a stronger equivalence relation than inner relation. Indeed, we can show that (lemma 4.11 in \cite{kellendonk}), if $\mathfrak{r}$ and $\mathfrak{s}$ are inner conjugate then it exists a unitary $w\in A$ such that they are also inner related with generator given by $w\mathfrak{r}(w^*)$.
An important result obtained in the paper \cite{kellendonk} is Corollary 4.14: it states that, under further further assumptions on the center of $A$ (satisfied in our case), given a reference real structure, up to stabilisation and inner conjugation, there is a finite number of real structures inner related to it. Moreover, if $\mathfrak{r}$ and $\mathfrak{s}$ are two inner conjugate real structures, we have that $(A, \gamma ,\mathfrak{r})$ and $(A, \gamma ,\mathfrak{s})$ are isomorphic as graded $C^*$-algebras equipped with a real structure: therefore, the $DK$-groups are isomorphic and we just have to focus on one real structure for each inner conjugation equivalence class. These are characterised by their relative signs.
\begin{definition}
\label{sign}
Let $(A,\gamma)$ be a graded $C^*$-algebra with two inner related real structures $\mathfrak{r}$, $\mathfrak{s}$. The
relative signs between $\mathfrak{r}$ and $\mathfrak{s}$ are
$$\eta_{\mathfrak{r}, \mathfrak{s}}=\left(\eta_{\mathfrak{r}, \mathfrak{s}}^1,\eta_{\mathfrak{r}, \mathfrak{s}}^2\right)\coloneqq\left(u\mathfrak{r}(u),u\gamma(u)^* \right)$$
where $u$ is any generator for $\mathfrak{s} \circ \mathfrak{r}$.
We say that a real symmetry $\mathfrak{r}$ is \textit{even} (\textit{odd}) if the relative sign $\eta_{\mathfrak{r}, \mathfrak{f}}^1$ to the reference structure $\mathfrak{f}$ is $+1$ ($-1$).
\end{definition}

In the case of chiral symmetry and TRS with real structure $\mathfrak{r}=\mathfrak{t}$, in the rough case we found that the $K$-group is $K_1(A^\mathfrak{t},\gamma)$. Now we can describe it more precisely. Indeed, we find that it is isomorphic to the $K$-groups of the following table. Moreover, if we assume the grading to be inner, we can get a more precise classification (see last column of the table).
\begin{equation}\label{chiralTRS_reference}
    \begin{array}{c c c |c c}
\hline K \text {-group } & \eta_{\mathfrak{r}, \mathfrak{f}} & K \text {-group }  & \text{grading} & K \text {-group (inner case) }  \\
&  & \text{w.r.t. reference}& & \text{w.r.t. reference} \\
\hline K_1(A^\mathfrak{t},\gamma) & (+1,+1)  & K_{1}\left(A^{\mathfrak{f}}, \gamma\right) & \text{real} & KO_1(A^{\mathfrak{f}}) \\K_1(A^\mathfrak{t},\gamma) & 
(+1,-1)  & K_{-1}\left(A^{\mathfrak{f} \circ \gamma}, \gamma\right) & \text{imag.} & KO_{-1}(A^{\mathfrak{f}}) \\K_1(A^\mathfrak{t},\gamma) & 
(-1,+1)  & K_{5}\left(A^{\mathfrak{f}}, \gamma\right) & \text{imag.} & KO_3(A^{\mathfrak{f}})\\K_1(A^\mathfrak{t},\gamma) & 
(-1,-1)  & K_{3}\left(A^{\mathfrak{f} \circ \gamma}, \gamma\right) & \text{real} & KO_5(A^{\mathfrak{f}})\\
\hline
\end{array}
\end{equation} 
\par
We consider now the case of a trivial grading with TRS or PHS and real structure $\mathfrak{r}$ or $\mathfrak{p}$ respectively. Again we tensor with $\C l_1$ to restore odd elements. The $K$-group is $KO_{0}\left(A^{\mathfrak{t}}\right) \text { or } K O_{2}\left(A^{\mathfrak{p}}\right) \text {, for TRS or PHS, respectively} $. We get 
\begin{equation}\label{TRS/PHS_ref}
    \begin{array}{ccccc}
\hline  \text { Symmetry } &K \text {-group } & \eta_{\mathfrak{r}, \mathfrak{f}}^{(1)}  & \mathfrak{s}  & K \text {-group } \\
& & &  & \text{w.r.t. reference}\\

\hline \text { TRS even } & KO_{0}\left(A^{\mathfrak{t}}\right)& +1 & +1& K O_{0}\left(A^{\mathfrak{f}}\right) \\
\text { TRS odd } & KO_{0}\left(A^{\mathfrak{t}}\right)&-1 & +1 & K O_{4}\left(A^{\mathfrak{f}}\right) \\
\text { PHS even } & K O_{2}\left(A^{\mathfrak{p}}\right) & +1  & -1 & K O_{2}\left(A^{\mathfrak{f}}\right) \\
\text { PHS odd } & K O_{2}\left(A^{\mathfrak{p}}\right) & -1  & -1  & K O_{6}\left(A^{\mathfrak{f}}\right) \\
\hline
\end{array}
\end{equation}
where we noted with $\mathfrak{s}=+1$ or $-1$ the real structure over $\C l_1 $ is $\mathfrak{l}_{1,0}$ or $\mathfrak{l}_{0,1}$ (the first case corresponds to a TRS symmetry, while the second to a PHS).

To sum up, when we assume the grading to be inner, we get the following classifications.
\begin{equation}\label{final_classifications}
\begin{array}{ccccc} 
\Xhline{5\arrayrulewidth} \text { Symmetries }  & \text { inner grading } & \eta_{\mathfrak{r}, \mathfrak{f}}^{(1)} & K \text {-group }& K \text {-group w.r.t. a reference}\\
& & & & \text {real structure }  \mathfrak{f}\\
\Xhline{5\arrayrulewidth}
\text { No Symmetry } & & & K U_{0}\left(A\right) & K U_{0}\left(A\right)\\
\hline
\text { Chiral } & & & K U_{1}\left(A\right) & K U_{1}\left(A\right)\\
\Xhline{3\arrayrulewidth}
\text { TRS } & & +1 &  KO_{0}\left(A^{\mathfrak{t}}\right) & KO_{0}\left(A^{\mathfrak{f}}\right)\\
 & & -1 &   & KO_{4}\left(A^{\mathfrak{f}}\right)\\
 \hline
\text { PHS } & & +1 &  KO_{2}\left(A^{\mathfrak{p}}\right) & KO_{2}\left(A^{\mathfrak{f}}\right)\\
 & & -1 &   & KO_{6}\left(A^{\mathfrak{f}}\right)\\
 \hline
\text { TRS + chiral }  & \text { real } & +1 & K O_{1}\left(A^{\mathfrak{t}}\right) & K O_{1}\left(A^{\mathfrak{f}}\right) \\
 &  & -1 &  & K O_{5}\left(A^{\mathfrak{f}}\right) \\
 \hline
 \text { TRS + chiral }& \text { imag. } & +1 & K O_{7}\left(A^{\mathfrak{t}}\right) & K O_{7}\left(A^{\mathfrak{f}}\right) \\
  &  & -1 &  & K O_{3}\left(A^{\mathfrak{f}}\right) \\
\hline
\end{array}
\end{equation}

\subsection{Two final remarks}
It is worth highlighting again two remarks. As already stated, the classifications presented above are relative and this is true from two distinct points of view. On the one hand, the Van Daele's construction depends on the choice of a reference element. This detail could seem innocuous as for two different choices of the reference element the groups which we obtain are isomorphic; however the isomorphism is not canonical. On the other hand, we remind that we can choose a reference real structure to get a finer classification and to be able to compare classifications with respect to other real structures related to the reference one (indeed this fact will be crucial to compare the classifications obtained with the two approaches). 

\section{Alldridge, Max and Zirnbauer's approach}
What follows relies on the paper \cite{zirnbauer} and the Ph.D. thesis of C. Max \cite{CMax} which develops the content of the paper.
We consider an approach to classify topological phases which gives a concrete and physically well justified $C^*$-algebra describing insulators/superconductors and symmetries. In particular, the power of the formalism presented in \cite{zirnbauer} is that we can consider in the same framework topological insulators (TI) and superconductors (TSC). The starting point is the description of our system (the insulator/superconductor) using the formalism of second quantization. The construction takes account of several approximations and hypothesis, namely we want to describe free fermions without (a priori) particle number conservation, using \textit{tight-binding approximation}, in a regime of \textit{homogeneous disorder} and allowing the possibility of having a homogeneous magnetic field.

The main building block of the construction is the so-called \textit{Nambu space} and (the Hamiltonian of) an insulator is a linear operator over this space. We will obtain then a Bogoliubov-de-Gennes (BdG) Hamiltonian: in this way, we have rewritten the description of the system in a first quantization language, which is standard and easier to deal with.
We call the relevant algebra (where the BdG Hamiltonians live) the \textit{algebra of bulk observables}. One of the achievement of \cite{zirnbauer} is to give a classification of topological phases (path-connected components of the set of BdG Hamiltonians possibly satisfying symmetry conditions). Also in this case Van Daele's version of $K$-theory is used.

To sum up, in this section we want to achieve three goals: the former is to build the $C^*$-algebra of \textit{bulk observables}, then to explain how to incorporate symmetries and finally showing which $K$-group are obtained for each case of symmetry.

\subsection{Fock and Nambu spaces}
We consider the Hilbert space $(\mathcal{H}, \langle \cdot,\cdot \rangle )$ which describes the behaviour of a single particle. We do not want to restrict ourselves to a single particle: we want to be able to consider an arbitrary number of particles at the same time. For this reason, we consider the tensor algebra of $\mathcal{H}$, where the tensor product with $x$ corresponds to adding a particle in the $x$ state to the system. Since we are working with fermions, Pauli's principle prescribes antisymmetrization and it is therefore natural to replace the tensor algebra with the exterior $\Lambda \mathcal{H}$ algebra (called the Fock space). We note the product on $\Lambda$ with $\wedge$. The exterior algebra is obtained as a tensor algebra modulo the ideal $I$ generated by the elements $v\otimes v$, with $v\in \mathcal{H}$.
The scalar product of $\mathcal{H}$ induces a scalar product $ \langle \cdot,\cdot \rangle_\Lambda $ on $\Lambda \mathcal{H}$. Indeed, for $n\in \N$ we set
$$\langle x_1\wedge x_2\wedge\dots\wedge x_n ,y_1\wedge y_2\wedge\dots\wedge y_n \rangle_\Lambda=\sum_{\sigma\in S_n}(-1)^{|\sigma|}\prod_{j=1}^{n} \langle x_j,y_{\sigma(j)} \rangle=\det(\langle x_i,y_j\rangle_{1\le i,j\le n})$$

and when $n\ne m$, $\langle x_1\wedge x_2\wedge\dots\wedge x_n ,y_1\wedge y_2\wedge\dots\wedge y_m \rangle_\Lambda=0$. Then we extend by linearity on the other elements.

To any $x\in \mathcal{H}$ we can associate a creation and annihilation operator $c^{\dagger}_x$ and $c_x$ such that

$$c^{\dagger}_x (x_1\wedge x_2\dots\wedge x_n)=x\wedge x_1\wedge x_2\dots\wedge x_n$$
and 
$$c_x(y_1\wedge y_2\dots\wedge y_n)=\sum_{i=1}^{n}(-1)^{i-1}\langle x,y_i\rangle y_1\wedge\dots\wedge y_{i-1}\wedge y_{i+1}\dots\wedge y_{n} $$
Alternatively, we can define $c^{\dagger}_x$ as the adjoint operator of $c_x$, i.e. the operator such that 
$$\langle c_x(y),w\rangle_\Lambda=\langle y,c^{\dagger}_x(w)\rangle_\Lambda\quad \text{for all }y\text{, }w \in \Lambda\mathcal{H}$$ 

We see that the application $x\mapsto c_x$ (resp. $x\mapsto c^{\dagger}_x$) from $\hil $ to $\mathcal{L}(\Lambda\mathcal{H})$ is anti-linear i.e. $c_{\text{i}x}=-\text{i}c_{x}$ (resp. linear). We note with $\rho:\mathcal{H}\rightarrow \mathcal{H}^*$ the Riesz isomorphism and we check easily that,thanks to antisymmetrization, for all $x\in \hil $ and $\varphi \in \hil^*$ we have
$$c_{x}^\dagger c_{x'}^\dagger=-c_{x'}^\dagger c_{x}^\dagger$$
$$
c_{\rho^{-1}(\varphi)} c_{\rho^{-1}(\varphi')}=-c_{\rho^{-1}(\varphi')} c_{\rho^{-1}(\varphi)}
$$
$$ c_{\rho^{-1}(\varphi')}c_{x}^\dagger=\varphi'(x)\text{Id}-c_{x}^\dagger c_{\rho^{-1}(\varphi')}$$

\begin{definition}[Field operator]
	A \textit{field operator} is an element of $\mathcal{L}(\Lambda\mathcal{H})$ (the linear operators on $\Lambda\mathcal{H}$) which is a $ \C$-linear combination of the operators $c^{\dagger}_x$ and $c_x$, with $x\in \mathcal{H}$.
\end{definition}
We note $\mathcal{W}\subset \mathcal{L}(\Lambda\mathcal{H})$ the complex vector space of \textit{field operators}. If we define the bilinear map $\eta$ as
$$\eta: \mathcal{H} \oplus \mathcal{H}^{*} \rightarrow \mathcal{W}, x+\varphi \mapsto c_{x}^{\dagger}+c_{\rho^{-1}(\varphi)}$$
we can see that $\mathcal{W}$ is isomorphic to $\mathcal{H}\oplus\mathcal{H}^*$ via $\eta$.
We define the anti-commutator $\{\cdot,\cdot \}:\mathcal{W}\times\mathcal{W}\rightarrow \mathcal{L}(\Lambda\mathcal{H})$
$$	\{\psi,\psi'\}\coloneqq\psi\psi'+\psi'\psi$$

\textbf{Remark:} $\mathcal{W}$ is not an algebra.

If $\psi=\eta(x+\varphi)$ and $\psi'=\eta(x'+\varphi')$ with $x, x'\in \mathcal{H}$ and $\varphi, \varphi'\in \mathcal{H}^*$, we can check that 
\begin{equation}
\tag{CAR relation}
\{\psi,\psi'\}=(\varphi(x')+\varphi'(x))\text{Id}=\tilde{q}(\psi,\psi')\text{Id}
\end{equation}
where we defined the bilinear form $\tilde{q}(\cdot,\cdot)=q(\eta^{-1}(\cdot),\eta^{-1}(\cdot))$ and $q(x+\varphi,x'+\varphi')=\varphi(x')+\varphi'(x)$. Indeed, we note that
$$\psi\psi'=(c_{x}^\dagger+c_{\rho^{-1}(\varphi)})(c_{x'}^\dagger+c_{\rho^{-1}(\varphi')})=c_{x}^\dagger c_{x'}^\dagger+c_{\rho^{-1}(\varphi)} c_{\rho^{-1}(\varphi')}+c_{x}^\dagger c_{\rho^{-1}(\varphi')}+c_{\rho^{-1}(\varphi)}c_{x'}^\dagger.$$

We have a natural real structure on $\mathcal{H}\oplus\mathcal{H}^*$:
\begin{equation}\label{NatRealStr}
    C=\begin{pmatrix}
0 & \rho^{-1} \\
\rho & 0
\end{pmatrix}
\end{equation}
and the corresponding structure over $\mathcal{W}$, $\gamma=\eta\circ C\circ \eta^{-1} $.
We can also define a scalar product over $\mathcal{H}^*$
$$\langle \varphi, \varphi' \rangle_{\mathcal{H}^*}=\langle \rho^{-1} (\varphi'),\rho^{-1} (\varphi) \rangle$$ 
The scalar product must be antilinear in the first argument and linear in the second. Moreover $\rho^{-1}$ is antilinear. We then exchange the order of $\varphi$ and $\varphi'$ to define a scalar product on $\mathcal{H}\oplus\mathcal{H}^*$ (and consequently on $\mathcal{W}$)

$$\langle x+\varphi, x'+\varphi' \rangle_{\mathcal{H}\oplus\mathcal{H}^*}=\langle x ,x' \rangle+\langle \varphi, \varphi' \rangle_{\mathcal{H}^*}$$
which defines then the structure of an Hilbert space on $\W$.

With this scalar product we can check that
\begin{equation}\label{eq:1}
	\langle \cdot,\cdot \rangle_{\mathcal{H}\oplus\mathcal{H}^*}=q(\gamma\cdot, \cdot)
\end{equation}

To summarise, we have three structures on $\mathcal{W}$
\begin{enumerate}
	\item the CAR relation (or $\tilde{q}$), coming from the canonical quadratic form;
	\item the real structure $\gamma$, coming from the Riesz isomorphism;
	\item the scalar product, coming from the scalar product on $\hil$.
\end{enumerate}
Previous structures are not independent, as shown in (\ref{eq:1}). With respect to these structures, one can define particular operators on $\mathcal{W}$:
\begin{enumerate}
	\item symmetric or anti-symmetric with respect to $\tilde{q}$ (i.e. $\tilde{q}(T\cdot,\cdot)=\pm\tilde{q}(\cdot,T\cdot)$), 
    \item real or imaginary with respect to $\text{Ad}_{\gamma}$ (i.e. $\text{Ad}_{\gamma}(T)= \pm T$).
	\item Self-adjoint or skew-adjoint with respect to $\langle \cdot, \cdot \rangle_{\mathcal{W}}$ ($\langle T\cdot, \cdot \rangle_{\mathcal{W}}=\pm\langle \cdot, T\cdot \rangle_{\mathcal{W}}$).
\end{enumerate}
In this case too, these definitions are not independent. 
\begin{definition}[Nambu space]
	The \textit{Nambu space} is the triple $\left(\mathcal{W}, \gamma 
 , q \right)$, where $q$ is the quadratic form giving rise to the CAR relation. 
\end{definition}
%

\subsection{Dynamics}
The evolution in time is described by a $\C$-linear $*$-automorphism $\mathcal{L}(\Lambda\mathcal{H})$. The dynamics is said to describe \textit{free fermions} if 
\begin{enumerate}
	\item the evolution preserves $\mathcal{W}$;
	\item evolution preserves the quadratic form (hence CAR), i.e.
	$$\tilde{q}(\psi(0),\psi'(0))=\tilde{q}(\psi(t),\psi'(t)), \forall t\ge 0$$
	\item evolution preserves the real structure, i.e. if $\psi(0)$ is real, then $\psi(t)$ is real for all $t\ge 0$.
\end{enumerate} 

We make the hypothesis that there exists $H\in\mathcal{L}(\Lambda\mathcal{H}) $ (called \textit{free fermion Hamiltonian}) such that if $A(0)\in\mathcal{L}(\Lambda\mathcal{H})$ then the time evolution of $A(t)$ is given by the Heisenberg equation
$$\dot{A}(t)=\text{i}\left[ H,A(t) \right]$$
This implies that 
$$A(t)=\text{Ad}_{\text{exp}(\text{i}tH)}(A(0))$$
To simplify the notation, we can define the Liouvillian $\mathfrak{L}\in \mathcal{L}(\mathcal{L}(\Lambda\mathcal{H}))$
$$\mathfrak{L}(\cdot)\coloneqq \text{i}[H,\cdot]$$
and the Heisenberg equation can be rewritten as
$$\dot{A}(t)=\mathfrak{L}(A(t))$$
$\mathfrak{L}$ is the generator of the time evolution.

The three previous conditions are equivalent to:
\begin{enumerate}
	\item $\mathfrak{L}(\mathcal{W})\subset\mathcal{W}$, which implies $[H,\mathcal{W}]\subset\mathcal{W}$.
	\item $\mathfrak{L}^T=-\mathfrak{L}$, i.e. $\tilde{q}(\mathfrak{L}\psi, \psi')=\tilde{q}(\psi, \mathfrak{L}\psi')$. 
	\item $\mathfrak{L}$ is real, i.e. $\text{Ad}_\gamma (\mathfrak{L})=\mathfrak{L}$. 
\end{enumerate}

\subsection{Bogoliubov-de-Gennes Hamiltonian}
Finally we define Bogoliubov-de-Gennes Hamiltonian (from now on, BdG Hamiltonian): as announced, we find a description of the physical system using the first quantization formalism, which is the standard framework to classify topological phases. 
The BdG Hamiltonian is defined as $\tilde{H}_{B}$
$$\tilde{H}_{B}\coloneqq -\text{i}\mathfrak{L}|_\mathcal{W}\in \mathcal{L}(\mathcal{W})$$
The BdG Hamiltonian can be related to its \textit{free fermion Hamiltonian}
$$\tilde{H}_{B}=[H,\cdot]$$
In \cite{zirnbauer} the BdG Hamiltonian is referred to as the Hamiltonian induced by the free fermion Hamiltonian. 
The Bogoliubov-de-Gennes Hamiltonian has the following properties:
\begin{enumerate}
	\item $\tilde{H}_{B}$ is self-adjoint with respect to the scalar product $\langle \cdot,\cdot\rangle_\mathcal{W}$. Indeed
	\begin{align*}
	    \langle \tilde{H}_{B}(\psi),\psi'\rangle_\mathcal{W}&=\langle -\text{i}\mathfrak{L}(\psi),\psi'\rangle_\mathcal{W}&\\
	    &=\text{i}\tilde{q}(\gamma\circ\mathfrak{L}(\psi),\psi')& \text{thanks to equation \ref{eq:1}}\\&=-\text{i}\tilde{q}(\mathfrak{L}^T(\gamma(\psi)),\psi')&\mathfrak{L}\text{ real and skew-symmetric}\\&=\langle \psi,\tilde{H}_{B}(\psi')\rangle_\mathcal{W}&
	\end{align*}
	
	\item $\tilde{H}_{B}$ is imaginary, i.e. $\text{Ad}_\gamma(\tilde{H}_{B})=-\tilde{H}_{B}$.
	Indeed
	$$\text{Ad}_\gamma(\tilde{H}_{B})=\gamma\circ(-\text{i}\mathfrak{L})\circ \gamma=+\text{i}\gamma\circ \mathfrak{L}\circ \gamma=+\text{i}\mathfrak{L}=-\tilde{H}_{B}$$
	($\gamma $ is anti-linear and $\mathfrak{L}$ is real).
\end{enumerate}
We can rewrite $\tilde{H}_{B}$ in matrix form with respect to $\mathcal{H}\oplus\mathcal{H}^*$, i.e. $H_{B}\coloneqq\eta^{-1}\circ\tilde{H}_{B}\circ\eta=\eta^{-1}([H,\eta(\cdot)])$ and $H_{B}$ has the form
$$H_{B}=\begin{pmatrix}
P & \Delta\\
\Delta' & P'
\end{pmatrix}$$
where $P\in \mathcal{L}(\mathcal{H})$, $\Delta\in \mathcal{L}(\mathcal{H}^*,\mathcal{H})$, $\Delta'\in \mathcal{L}(\mathcal{H},\mathcal{H}^*)$ and $P'\in \mathcal{L}(\mathcal{H}^*)$.

It follows from the above properties that:
\begin{itemize}
	\item $H_{B}$ is self-adjoint $\iff $ $(\Delta^*=\Delta')\text{ and } \left(P \text{ and } P' \text{ are self-adjoint}\right)$.
	\item $H_{B}$ is imaginary w.r.t. the natural real structure (i.e. $CH_{B}C^{-1}=-H_{B}$, see (\ref{NatRealStr})) $\iff$ $\Delta'=-\rho \Delta \rho\text{ and } P'=-\rho P\rho^{-1}$.
\end{itemize}
Therefore 
$$H_{B}=\begin{pmatrix}
P & \Delta\\
-\rho \Delta \rho & -\rho P\rho^{-1}
\end{pmatrix}$$
with $P$ self-adjoint and $\Delta^*=-\rho \Delta \rho$.
A BdG Hamiltonian $H_{B}$ is \textit{gapped} if it is invertible. A gapped BdG Hamiltonian is flattened (and so unitary, as $H_{B}$ is self-adjoint) if $H_{B}^2=1$. Moreover, we can associate to any gapped BdG Hamiltonian $H_{B}$ a flattened BdG Hamiltonian $H_{B}\lvert H_{B}\rvert ^{-1}$.

The power of the Nambu space formalism is that we can consider in the same framework  topological insulators (TI) and superconductors (TSC): we can identify TI and TSC with their flattened BdG Hamiltonian. For TI it easy to understand the requirement of being gapped. Indeed, the case of TI corresponds to the situation where we have conservation of the charge: we will see that in this case we can reduce the algebra describing multi-particles and we will recover the single-particle setting. In this framework, the physical meaning of a gapped Hamiltonian is well-motivated. Regarding TSC, even if it is not completely clear to us how to justify the existence of a spectral gap in the general case, in several specific cases one can relate a gap in the spectrum of $H_B$ to a superconducting gap.

\begin{definition}[Abstract BdG Hamiltonian]
	An \textit{abstract BdG Hamiltonian} is a self-adjoint, skew-symmetric, invertible and unitary element $H_{B}\in\mathcal{L}(\W)$ which is imaginary with respect to the canonical real structure $\text{Ad}_\gamma$.
\end{definition}

\textbf{Remark: comparison with the BdG Hamiltonian coming from physics.}\\
We consider for simplicity the case where $\hil$ is a Hilbert space of dimension $n$ and we note with $\mathbf{c}^\dagger=\left( c_{1}^\dagger,c_{2}^\dagger, \cdots c_{n}^\dagger \right) $ (resp. $ \mathbf{c}$) the vectors containing the creation (resp. annihilation) operators in the states corresponding to the element of a canonical basis of $\hil$. In physics literature (\cite{el-batanouny_2020}), we find that $H$ has the general form $$H=\frac{1}{2}(\mathbf{c}^\dagger, \mathbf{c})\mathfrak{H}(\mathbf{c}, \mathbf{c}^\dagger)$$
where
$$\mathfrak{H}=\begin{pmatrix}
\Theta & \Xi \\
\tilde{\Xi} & \Gamma
\end{pmatrix}$$
Explicitly 
\begin{equation}\label{eq2:BdG}
    H=\frac{1}{2}\left( \sum_{l,m=1,\cdots,n}
    c_l^\dagger\Theta_{l,m}c_m+c_l^\dagger\Xi_{l,m}c_m^\dagger + c_l\tilde{\Xi}_{l,m}c_m+c_l\Gamma_{l,m}c_m^\dagger\right) 
\end{equation}
We can see that the Hamiltonian BdG corresponding to $H$ is $\mathfrak{H}$. Therefore we must have
$$\mathfrak{H}=\begin{pmatrix}
\Theta & \Xi \\
-\rho\Xi\rho^{-1} & -\rho \Theta\rho^{-1}
\end{pmatrix}$$
with $\Theta$ self-adjoint and $\Xi^*=-\rho \Xi \rho$.

\subsection{The $C^*$-algebra of disordered bulk observables}\label{BulkAlgebra}

Until now, we did not specify the single-particle complex Hilbert space $\hil$. We describe a solid state system as a $d$-dimensional crystalline lattice in $\R^d$ and taking account of possible internal degrees of freedom (including the description of a basis of the unit cell, spin, orbital degrees of freedom...). Therefore it is reasonable to define
$$\hil\coloneqq \ell^2(\Lambda)\otimes V$$
where $\Lambda$ is the Bravais lattice (isomorphic to $\R^d$) and $V$ is a finite-dimensional complex Hilbert space. If we define $W\coloneqq  V\oplus V^*$ (therefore $W$ is the Nambu space corresponding to the local finite-dimensional Hilbert space $V$) we can show the following result.
\begin{theorem}\label{lemmaRealStructure1}
    The $C^{*,r}$-algebra $(\mathcal{L}(\W), \text{Ad}_\gamma$ is $*$-isomorphic to $\mathcal{L}(\ell^2(\Lambda))\otimes \mathcal{L}(W)$ with real structure $\overline{T_1\otimes T_2}=\mathfrak{c}_\Lambda T_1 \mathfrak{c}_\Lambda\otimes \gamma T_2 \gamma$ ($T_1\in \mathcal{L}(\ell^2(\Lambda))$ and $T_2\in \mathcal{L}(W)$) where $\mathfrak{c}_\Lambda$ denotes the point-wise complex conjugation on $\ell^2(\Lambda)$.
\end{theorem} 
In the following, we will often talk about \textit{quaternionic structures} over a Hilbert space.
\begin{definition}
An anti-linear unitary operator over a Hilbert space $\hil$ $T : \hil \rightarrow \hil$ is called a \textit{quaternionic structure} if $T^2=-1$. That is, T is
anti-linear and we have
$$\langle T x, T y\rangle=\langle y, x\rangle, \quad\langle T x, y\rangle=-\langle T y, x\rangle \quad \forall x, y \in \mathcal{H}$$
A complex Hilbert space equipped with a quaternionic structure is called a quaternionic Hilbert space. A linear map between quaternionic Hilbert spaces is called quaternionic if it intertwines the quaternionic structures.
\end{definition}

In the same way of the lemma, if $W$ has a real or quaternionic structure $\mathfrak{r}$ (commuting with $\gamma$), we can consider an \textit{induced} real structure on $\mathcal{L}(\ell^2(\Lambda))\otimes \mathcal{L}(W)$ as $\overline{T_1\otimes T_2}=\mathfrak{c}_\Lambda T_1 \mathfrak{c}_\Lambda\otimes \mathfrak{r} T_2 \mathfrak{r}$ ($T_1\in \mathcal{L}(\ell^2(\Lambda))$ and $T_2\in \mathcal{L}(W)$).

In the framework of the tight-binding approximation, we are interested in Hamiltonians with a finite hopping range, namely to \textit{controlled operators}.
\begin{definition}[Controlled operator]
    We call an operator $O\in \mathcal{L}(\ell^2(\Lambda))\otimes \mathcal{L}(W)$ a \textit{controlled operator} if it exists $R>0$ such that for all $x$, $y\in \Lambda$ with $\lVert x-y \rVert >R$ we have
    $$\langle \delta_x\otimes w, O(\delta_y\otimes w') \rangle=0 \quad \forall w, w'\in \W $$
    where $\delta_x$ denotes the element of $\ell^2(\Lambda)$ which takes the value $1$ on the site $x$ and $0$ elsewhere.

    We denote with $C_u^*(\Lambda, W)$ the norm closure of the set of controlled operators in $\mathcal{L}(\ell^2(\Lambda))\otimes \mathcal{L}(W)$, which is a $C^*$-algebra.
\end{definition}
$C_u^*(\Lambda, W)$ is also called the \textit{uniform Roe algebra}.
The following lemma allows us to extend naturally real and quaternionic structures from $W$ to $\mathcal{L}(\ell^2(\Lambda))\otimes \mathcal{L}(W)$.
\begin{theorem}\label{lemmaRealStructure2}
    $C_u^*(\Lambda, W)$ is a $C^*$-algebra and if $W$ has a real or quaternionic structure $C_u^*(\Lambda, W)$ is a $C^{*,r}$-algebra with respect to the induced real structure on $\mathcal{L}(\ell^2(\Lambda))\otimes \mathcal{L}(W)$.
\end{theorem}

Now we consider homogeneous disorder. Asking for a solid state system to be microscopically translation invariant is not a reasonable assumption, as all such systems display some disorder which breaks the translational invariance. One may however assume that the translational invariance still holds at a macroscopic level. We introduce the space of disorder configurations.
\begin{definition}[Space of disorder configurations]
	A space of disorder configurations is a compact Hausdorff space $\Omega$ with a continuous action of the Abelian group $\Z^d$. We denote the action by the right action $\omega \mapsto \omega \cdot x$ for $\omega \in \Omega$ and $x \in \Z^d$. Furthermore, we assume that there is a fully supported $\Z^d$-invariant Borel probability measure $\mathbb{P}$ on $\Omega$.
\end{definition}
If the Hamiltonian is macroscopically homogeneous, then, together with its translates, it should be a function on $\Omega$ which is translationally covariant in the sense that it intertwines the action
of the translation group on $\Omega$ and the space of operators.

Finally, we deal with homogeneous magnetic fields. Magnetic field is modelled by a real skew-symmetric $d\times d$-matrix $\mathbf{B} = (B_{\mu\nu})$. In our context we define 
$$\sigma(x,y)\coloneqq \text{e}^{\frac{\text{i}}{2}(x, \textbf{B}y)}\quad \forall x,y \in \Lambda$$
and its corresponding \textit{magnetic translations} $u_x^\sigma$ on $\hil$
$$\left(u_x^\sigma \psi\right)(y):=\sigma(y, x) \psi(y-x), \quad \forall y \in \Lambda, \psi \in \mathcal{H}$$
for $x\in \Lambda$ with integer coefficients with respect to the vectors defining the lattice. By the canonical lift, $\sigma$ can be extended to $\W$ as follows
$$
\tilde{\sigma}(x, y):=\left(\begin{array}{cc}
\text{e}^{\frac{i}{2}(x, \mathbf{B} y)} & 0 \\
0 & \text{e}^{-\frac{i}{2}(x, \mathbf{B} y)}
\end{array}\right)
$$
and the extension is real with respect to $\text{Ad}_\gamma$.
Once we have defined controlled operators, homogeneous disorder and the magnetic field, we can define the algebra containing physically relevant operators.
\begin{definition}[$C^*$-algebra of bulk observables]\label{defBulkAlg}
We define $A^{W,\tilde{\sigma}}$ to be the set of all maps $O : \Omega \rightarrow C_u^*(\Lambda, W)$ that are continuous in the norm topology and covariant in the sense that 
$$O_{\omega \cdot x}=(u_x^{\tilde{\sigma}})^* O_\omega u_x^{\tilde{\sigma}}, \quad \forall x \in \Z^d, \omega\in\Omega$$
With point-wise operations 
\begin{align*}
    \left(O_1+\lambda O_2\right)_\omega & :=\left(O_1\right)_\omega+\lambda\left(O_2\right)_\omega, \quad , \forall O, O_1, O_2 \in A^{W, \sigma}, \lambda \in \mathbb{C}, \omega \in \Omega\\ 
    \left(O_1 O_2\right)_\omega & :=\left(O_1\right)_\omega \left(O_2\right)_\omega\\ \left(O^*\right)_\omega & :=\left(O_\omega\right)^*,
\end{align*}
and norm defined by 
$$\lVert O \rVert \coloneqq \text{sup}_{\omega\in\Omega} \lVert O_\omega \rVert, \quad \forall O\in A^{W,\tilde{\sigma}},$$
the set $A^{W,\tilde{\sigma}}$ is a $C^*$-algebra. If $W$ is real or quaternionic and $\tilde{\sigma}\equiv 1$, then $A^W$ is a real $C^*$-algebra with real structure defined by
$$(\overline{O})_\omega\coloneqq \overline{(O_\omega)}, \quad \forall O\in A^W, \omega \in \Omega$$
The algebra of \textit{bulk observables} is defined as the real $C^*$-algebra
$$\mathbb{A}\coloneqq A^W.$$
In case of charge-conserving observable, the algebra $A^{W,\tilde{\sigma}}$ can be reduced (we will see later why) and therefore it is useful to define the algebra of complex bulk observables as the complex $C^*$-algebra
$$\text{A}\coloneqq A^{V,\sigma}.$$
\end{definition}
\textbf{Remark.} In the real case we do not consider any non-trivial magnetic fields, as an homogeneous magnetic field is not compatible with a real symmetry in general. On the other hand, we can do it in the case of charge conservation (i.e. what we call the complex case). The reason is that, even if we begin with an algebra with a real structure (the natural real structure $\gamma$), this particular real structure is compatible with a non-trivial magnetic field (\textit{cf.} discussion in \cite{zirnbauer}). 

\subsection{Symmetries}\label{subsec:symmetries}
In this approach, symmetries are local, namely we just consider symmetries acting on the local part of the Nambu space $\W$ (namely on $W=V\oplus V^*$ and not $\ell^2(\Lambda)$) through linear or anti-linear real unitary isomorphisms. Lemmas \ref{lemmaRealStructure1} and \ref{lemmaRealStructure2} show how to obtain the induced real structure on $\mathcal{L}(\ell^2(\Lambda))\otimes \mathcal{L}(W)$ from a real or quaternionic structure on $W$. Then, it can be extended (following the definition of bulk algebra) to $\mathbb{A}$. Linear isomorphisms on $W$ are trivially extended to $\mathbb{A}$ and A (the isomorphism is not trivial just on the $W$-part). For these reasons and in order to have a lighter notation, in the following we just consider symmetries on $W$.
\begin{definition}[Symmetry]
	A linear or anti-linear real unitary isomorphism $S : V\oplus V^* \rightarrow V\oplus V^*$ is called a \textit{symmetry} and the BdG Hamiltonian $H_{B}$ is called $S$-symmetric if
    $$[H_{B}, S]=0$$
\end{definition}
We consider four examples of symmetries: with suitable combinations of these symmetries one can represent all of the ten symmetry classes in the Tenfold Way \cite{zirnbauer}.
\begin{enumerate}
        \item \textbf{Time Reversal Symmetry} (TRS). It is defined by a anti-linear unitary automorphism $T:V\rightarrow V$ such that $T^2=\pm 1$ (even or odd symmetry). In the case of fermions, we have $T^2=-1$ (i.e. $T$ defines a quaternionic structure over $V$). This is the reason why we will always (when there is no Spin Rotations Symmetry) observe an odd TRS in the following (\textit{cf.} first lines of the final table). By the canonical lift, we extend $T$ to $W$ with
        $$\Tilde{T}\coloneqq \begin{pmatrix}
T & 0 \\
0 & \rho T \rho^{-1}
\end{pmatrix}$$
We easily check that the operator $\Tilde{T}$ is real with respect to the natural real structure $\text{Ad}_\gamma$ ($\text{Ad}_\gamma(\Tilde{T})=\Tilde{T}$). Moreover, the operator $\gamma T$ is skew-hermitian, indeed
\begin{align*}
  \langle (\gamma T)^*(v_1,\phi_1),(v_2,\phi_2)\rangle&=\langle (v_1,\phi_1), (\gamma T)(v_2,\phi_2)\rangle \\
  &= \langle T(v_2),\rho^{-1}\phi_1\rangle+\langle v_1,T\rho^{-1}\phi_2\rangle\\
  &=-\langle T\rho^{-1}\phi_1,v_2\rangle-\langle \rho^{-1}\phi_2,Tv_1\rangle\\&=-\langle (\gamma T)(v_1,\phi_1),(v_2,\phi_2)\rangle
\end{align*}
\item \textbf{Spin rotations symmetries} (SRS). We are dealing with electrons, therefore with particles with $\frac{1}{2}$-spin. $\frac{1}{2}$-spin is described as the lowest irreducible unitary representation of the Lie group $SU(2)$ over the single-particle space $V$. Spin Rotations Symmetry means to impose the commutation of the Hamiltonian with a set of generators of a representation of the Lie algebra $\mathfrak{su}(2)$ on $V$.
We define the generators of a representation on $V$ as $j_1=\text{i}S_1$, $j_2=\text{i}S_2$ and $j_3=\text{i}S_3$, such that they satisfy the commutation relations
$$[j_2, j_1]=2j_3 \quad [j_1, j_3]=2j_2 \quad [j_3, j_2]=2j_1 $$ 
and we assume moreover that $S_i^2$'s square to $1$. This assumption is reasonable because the spin rotations operators act on the fiber $\C^2$ of $V=V'\otimes \C^2$. As the representation must be compatible with the $*$-operation, we have that $j_i$'s must be skew-adjoint (and therefore the $S_i$'s self-adjoint). 
Since time-reversal $T$ inverts the spin we also assume that $T$ anti-commutes with the $S_1$, $S_2$, $S_3$ and subsequently it commutes with $j_1$, $j_2$, $j_3$.

The canonical lift to symmetries on $W$ is given by
$$\Tilde{j_\mu}\coloneqq \begin{pmatrix}
j_\mu & 0 \\
0 & \rho j_\mu \rho^{-1}
\end{pmatrix}\in \mathcal{L}(W)$$
We easily check that the operators $\Tilde{j}_\mu$ are real with respect to the natural real structure $\text{Ad}_\gamma$.

\item \textbf{Charge conservation} (Q). The charge operator is defined as
$$Q=\begin{pmatrix}
1 & 0 \\
0 & -1
\end{pmatrix}\in \mathcal{L}(W)$$
Indeed, 
$$H_{B}=\begin{pmatrix}
P & \Delta\\
\Delta^* & -P^T
\end{pmatrix}$$
preserves the charge if and only if $\Delta=0$ (the off-diagonal blocks correspond to the terms increasing and decreasing of two units the charges, see eq. \ref{eq2:BdG}). Therefore, the operator $Q$ is a diagonal block matrix. Moreover the upper block and the lower one must have opposite signs as the former corresponds to particles and the latter to anti-particles.
\item \textbf{Particle-hole symmetry} (PHS). Consider a unitary operator $S:V\rightarrow V$ such that $S^2=1$, $[S,T]=0$ and $[S,S_\mu]=0$. We define the symmetry
$$L=\begin{pmatrix}
0 & S\rho^{-1} \\
\rho S & 0
\end{pmatrix}$$
$L$ is anti-linear. We observe that $L$ commutes with $j_\mu$ (consequence of $[S,S_\mu]=0$). The name PHS is justified as we see that $L$ is off-diagonal, therefore it changes a particle in a anti-particle (a hole) and vice-versa.
    \end{enumerate}
\subsection{Tenfold Way and $K$-groups}
Once the framework is set, we want to describe topological phases associated to each case of symmetry. Namely, we want to say that two BdG Hamiltonians are equivalent (they are in the same topological phase) when we can find a continuous path of invertible, imaginary elements connecting the two Hamiltonians which moreover preserves symmetries (if any). 
As we saw for the abstract approach, the tool to obtain this result is Van Daele's version of $K$-theory: starting from the real and complex bulk algebras $\mathbb{A}$ and A and imposing some symmetries, we obtain a $K$-group where each element corresponds to a topological phase. With these definitions, in \cite{zirnbauer} the authors obtain the following table:
$$\begin{array}{cccc}
\hline s & \text { Class } & \text { symmetries } & K-\text {group}\\
\hline 
0 & \text { D } & \text {none} & KO_2(\mathbb{A}^{\text{Ad}_\gamma})\\
1 & \text { DIII } & \text {TRS} & KO_{3}(\mathbb{A}^{\text{Ad}_\gamma})\\
2 & \text { AII } & \text {TRS+Q} & KO_{4}(\mathbb{A}^{\text{Ad}_\gamma})\\
3 & \text { CII } & \text {TRS+Q+PHS} & KO_{5}(\mathbb{A}^{\text{Ad}_\gamma})\\
4 & \text { C } & \text {SRS} & KO_{6}(\mathbb{A}^{\text{Ad}_\gamma})\\
5 & \text { CI } & \text {SRS+TRS} & KO_{7}(\mathbb{A}^{\text{Ad}_\gamma})\\
6 & \text { AI } & \text {SRS+TRS+Q} & KO_{0}(\mathbb{A}^{\text{Ad}_\gamma})\\
7 & \text { BDI } & \text {SRS+TRS+Q+PHS} & KO_{1}(\mathbb{A}^{\text{Ad}_\gamma})\\
\hline
0 & \text{A} & \text { Q } & KU_0(\text{A})\\
1 & \text { AIII } & \text { Q+PHS } & KU_1(\text{A})\\
\end{array}
$$

The first eight lines correspond to the real cases, while the last two lines correspond to the complex cases. We explain the meaning of the column $s$. In \cite{zirnbauer}, for technical reasons (namely to deal in a compact way with the fact that operators describing symmetries can be linear or anti-linear), the authors introduce the definition of \textit{pseudo-symmetry} and \textit{quasi-paricle vacua} (essentially the BdG Hamitonian multiplied by i). The pseudo-symmetries are obtained from symmetries defined before and the column $s$ corresponds to the number of pseudo-symmetries they use. Concerning the column "Class", it corresponds to the historical names of topological phases (related to the classification of symmetric spaces) as presented in the Kitaev table \cite{kitaev2009periodic}.

\section{Comparison of the two approaches}

We would like to see how the two approaches \cite{kellendonk} and \cite{zirnbauer} are related. Our goal is to see \cite{zirnbauer} as a particular case (a specialisation) of the framework of \cite{kellendonk}. In particular,\cite{zirnbauer} proposes ten combinations of symmetries (two complex cases and eight real cases) and for each case they computes a $K$-group. We show that we can translate the relations in terms of chiral symmetry, TRS, PHS as in \cite{kellendonk} and that we find the same $K$-groups. In particular in this way we can bypass the technical construction of \cite{zirnbauer} of Van Daele's $K$-groups using pseudo-symmetries and quasi-particle vacua and therefore, we give a clearer view of the concrete approach. 
\subsection{Preliminaries}
As explained in paragraph \ref{subsec:symmetries}, the relevant algebras for us are the real and complex (charge conserving) bulk algebras $\mathbb{A}$ and A. Once again we remind that symmetries concern just the internal degrees of freedom: then, starting from a symmetry on $W$, we can naturally induce the corresponding symmetry on $\mathbb{A}$ and A. For this reason, in order to make notation lighter, we can consider as relevant algebras 
$$\mathbb{A}_{\text{loc}}=\mathcal{L}(W) \quad\text{ and }\quad \text{A}_{\text{loc}}=\mathcal{L}(V)$$
where "loc" stands for local (i.e. we jsut consider internal degrees of freedom).
At the end, we will be able to replace $\mathbb{A}_{\text{loc}}$ (with real structure $\text{Ad}_\gamma$) and $\text{A}_{\text{loc}}$ with $\mathbb{A}$ and A without further concerns. 

We will often use the isomorphism
\begin{equation}\label{M2(A)}
    \psi:\mathbb{A}_{\text{loc}}\xrightarrow{\sim} M_2(\text{A}_{\text{loc}})
\end{equation}
as $W\cong V\oplus V^*$. The isomorphism is given by
\begin{align*}
  \psi\colon M_2(\mathcal{L}(V)) & \longrightarrow \mathcal{L}(W) \\
  \begin{pmatrix}
\phi_1 & \phi_2 \\
\phi_3 & \phi_4 
\end{pmatrix} & \longmapsto
  \begin{pmatrix}
\phi_1 & \phi_2 \circ\mathfrak{m} ^{-1}\\
\mathfrak{m} \circ\phi_3 & \mathfrak{m} \circ\phi_4 \circ\mathfrak{m} ^{-1}
\end{pmatrix} 
\end{align*}
where $\mathfrak{m}\colon V  \longrightarrow V^*$ is defined as the Riesz isomorphism on the elements of a basis of $V$ and then extended $\C$-linearly on all $V$. Therefore, the map $\mathfrak{m}$ depends on the choice of a basis for $V$, but this is not a problem.

\textbf{Remark.} We want $\psi$ to be a linear isomorphism: $\mathfrak{m}$ is not the Riesz isomorphism, which is anti-linear!

We verify easily that 
$$\psi(\cdot)=\text{Ad}_{\mathfrak{M}}(\cdot)\quad \text{with}\quad \mathfrak{M}=\begin{pmatrix}
    1 & 0\\
    0 & \mathfrak{m}
\end{pmatrix}$$

This isomorphism will be useful in the following in order to compare the $K$-groups of the real algebras $\text{A}_{\text{loc}}^\mathfrak{r}$ and $\mathbb{A}_{\text{loc}}^{\Tilde{\mathfrak{r}}}$, where $\mathfrak{r}$ is a real structure and $\Tilde{\mathfrak{r}}$ is the natural extension to $\mathbb{A}_{\text{loc}}$
$$\Tilde{\mathfrak{r}}=\begin{pmatrix}
    \mathfrak{r} & 0\\
    0 & \rho \mathfrak{r} {\rho}^{-1}
\end{pmatrix}$$
In particular, we will often use the following result:

\begin{theorem}\label{lemma1}
    Given the time reversal operator $T$ over $V$ and the corresponding operator $\Tilde{T}$ over $ W$, we have an isomorphism between the $i$-th real $K$-groups
    
    $$KO_i(\text{A}_{\text{loc}}^{\text{Ad}_T})\cong KO_i ( \mathbb{A}_{\text{loc}}^{\text{Ad}_{\Tilde{T}}} )$$
\end{theorem}
\begin{proof}
    We set $\mathfrak{r}=\text{Ad}_T$ and $\tilde{\mathfrak{r}}=\text{Ad}_{\Tilde{T}}$. Moreover, we define $\mathfrak{r}_2$ as the real structure over $M_2(\text{A}_{\text{loc}})$ where we apply $\mathfrak{r}$ element-wise. The real algebras $\text{A}_{\text{loc}}^{\mathfrak{r}}$ and $M_2(\text{A}_{\text{loc}})^{{\mathfrak{r}}_2}$ are Morita equivalent and therefore, thanks to the stability of the $K$-groups with respect to Morita equivalence (as reminded in \cite{kellendonk}, page 2282), we have
    $$KO_i(\text{A}_{\text{loc}}^{\mathfrak{r}})\cong KO_i ( M_2(\text{A}_{\text{loc}})^{\mathfrak{r}_2})$$
    We just need to show that $M_2(\text{A}_{\text{loc}})^{\mathfrak{r}_2}$ and $\mathbb{A}_{\text{loc}}^{\Tilde{\mathfrak{r}}}$ are isomorphic. The isomorphism of real algebras is given by $\psi$ defined above, because $\psi$ intertwines the real structures, i.e.
    $\psi({\mathfrak{r}_2}(M)) =\Tilde{\mathfrak{r}}(\psi(M))$. To show it, first we observe that 
    $$\mathfrak{r}_2=\text{Ad}_{\mathfrak{R}}\quad \text{with}\quad \mathfrak{R}=\begin{pmatrix}
    T & 0\\
    0 & T
\end{pmatrix}$$
and that 
$$\Tilde{T}=\begin{pmatrix}
    T & 0\\
    0 & \rho T{\rho}^{-1}
    \end{pmatrix}=\begin{pmatrix}
    T & 0\\
    0 & \mathfrak{m}T{\mathfrak{m}}^{-1}
\end{pmatrix}$$
because of the definition of $\mathfrak{m}$.
Then, for all $M\in M_2(\text{A}_{\text{loc}})$     \begin{equation*}
        \begin{split}
\psi({\mathfrak{r}_2}(M)) & = \begin{pmatrix}
    1 & 0\\
    0 & \mathfrak{m}
\end{pmatrix}\begin{pmatrix}
    T & 0\\
    0 & T
\end{pmatrix} M \begin{pmatrix}
    T^{-1} & 0\\
    0 & T^{-1}
\end{pmatrix}\begin{pmatrix}
    1 & 0\\
    0 & \mathfrak{m}^{-1}
\end{pmatrix}\\
 & = \begin{pmatrix}
    T & 0\\
    0 & \mathfrak{m}T\mathfrak{m}^{-1}
\end{pmatrix}\begin{pmatrix}
    1 & 0\\
    0 & \mathfrak{m}
\end{pmatrix} M \begin{pmatrix}
    1 & 0\\
    0 & \mathfrak{m}^{-1}
\end{pmatrix}\begin{pmatrix}
    T & 0\\
    0 & \mathfrak{m}T\mathfrak{m}^{-1}
\end{pmatrix}^{-1}\\&=\Tilde{\mathfrak{r}}(\psi(M))
\end{split}
    \end{equation*}
\end{proof}

The physical objects we are interested in are the Hamiltonians of Bogoliubov-de-Gennes (BdG), which are described mathematically as self-adjoint, imaginary (with respect to the natural real structure $\gamma$), unitary elements of $\mathbb{A}_{\text{loc}}$ satisfying a set of symmetry conditions. 
We introduce here the notation to take into account symmetries. Given an order two $*$-isomorphism $r$ which could be linear or anti-linear, we say that a BdG Hamiltonian $H_B$ satisfies the relation defined by $r$ if $r(H_B)=H_B$. If we define the set of relations $\mathcal{R}=\{r_1, \dots, r_k\}$, we say that $H_B$ satisfies the relations $\mathcal{R}$ if $H_B$ satisfies the relation defined by $r_i$ for $ i=1,\dots k$. We take care of the condition of being imaginary (w.r.t. the natural real structure) with one of the $r_i$'s.

We note the set of BdG Hamiltonians up to homotopy (therefore the set of topological phases) by
$$GL^{s.a.}_{h}(\mathbb{A}_{\text{loc}}, \mathcal{R})$$
as our Hamiltonians are invertible and self-adjoint elements.
As announced in the introduction, our goal is to describe $GL^{s.a.}_{h}(\mathbb{A}_{\text{loc}}, \mathcal{R})$. To do that, we assign to it a Van-Daele's $K$-group
$$GL^{s.a.}_{h}(\mathbb{A}_{\text{loc}}, \mathcal{R})\leadsto DK(\mathbb{A}_{\text{loc}}', \alpha, \mathfrak{r})$$
where $\alpha$ and $\mathfrak{r}$ are a suitable grading and real structure respectively and where $\mathbb{A}_{\text{loc}}'$ is a new algebra to take care of all relations of $\mathcal{R}$ (in our case, $\mathbb{A}_{\text{loc}}'$ will be the tensor product between $\mathbb{A}_{\text{loc}}$ and a suitable Clifford algebra). The goal of our work is to find for each case of symmetry the appropriate choice for $\mathbb{A}_{\text{loc}}'$, $\alpha$ and $\mathfrak{r}$: in what follows we will show explicit calculations to find them.
An important question is what the meaning of this association is, namely what the "$\leadsto$" symbol stands for. We saw in the construction of Van-Daele's $K$-groups that the group $DK(A,\alpha)$ corresponds quite well to the set of path-connected components of odd, self-adjoint, invertible elements of $A$, but with an important difference: the classification is up to stabilisation! Stabilisation means that we consider homotopies in the matrix algebra of $A$, where the matrix algebra is the set of matrices of arbitrary size with entries in $A$. Therefore, $GL^{s.a.}_{h}(A, \{-\alpha\})$ and $DK(A,\alpha)$ are not the same object, but we can justify the association. From the physics point of view, the fact of considering the "stabilisation" corresponds to the stacking of the insulator (multiple layers of insulator) and therefore it is physically reasonable to consider the matrix algebra. On the other hand, from the mathematics perspective, algebraic topology gives us several tools to compute more easily the group $ DK(A, \alpha)$ (rather than $GL^{s.a.}_{h}(A, \{-\alpha\})$).

We divide the symmetry cases in three categories: complex cases, real cases with spin rotations symmetry and real cases without spin rotations symmetry. The complex cases have in common the condition of charge conservation which allows us to reduce the algebra and just consider the reduced complex algebra $\text{A}_{\text{loc}}$. As a final remark, we anticipate that the symmetry with respect to the operators associated to spin rotations plays a special role: for this reason, we treat the real cases with and without SRS in separate paragraphs.

\subsection{Complex cases}
A BdG Hamiltonian has the general form
$$H_{B}=\begin{pmatrix}
P & \Delta\\
-\rho \Delta \rho & -\rho P\rho^{-1}
\end{pmatrix}\in \mathcal{L}(V\oplus V^*)$$
with $\Delta$ anti-symmetric and $P$ self-adjoint. 
If $H$ commutes with $Q$, we have
$$H_{B}=\begin{pmatrix}
P & 0\\
0 & -\rho P\rho^{-1}
\end{pmatrix}$$
The condition of being imaginary is already satisfied: we just have two copies of the same operator. Therefore, we can restrict ourselves to the reduced complex algebra $\text{A}_{\text{loc}}$ (complex because there is no more a real structure). We consider here the two symmetry cases which allows to deal with a complex $C^*$-algebra. 

\subsubsection{Charge conservation}
$$\boxed{\mathcal{R}=\{-\text{Ad}_\gamma,\text{Ad}_Q\}}$$
As explained just above, we consider the reduced complex algebra $\text{A}_{\text{loc}}$ and the real condition is already satisfied. We are in the case of no symmetry of table \ref{final_classifications}, therefore $$GL^{s.a.}_{h}(\mathbb{A}_{\text{loc}}, \mathcal{R})=GL^{s.a.}_{h}(\text{A}_{\text{loc}}) \leadsto DK(\text{A}_{\text{loc}}\otimes\C l _1, \text{id}\otimes\text{st})\cong KU_0(\text{A}_{\text{loc}})$$

\subsubsection{Charge conservation + particle hole symmetry}
$$\boxed{\mathcal{R}=\{-\text{Ad}_\gamma,\text{Ad}_Q, \text{Ad}_L \}}$$

where 
$$L=\begin{pmatrix}
0 & S\rho^{-1} \\
\rho S & 0
\end{pmatrix}$$
and $S\in \mathcal{L}( V)$ a unitary operator.
As before, charge conservation and the condition of being imaginary allow us to consider just the block $P\in \mathcal{L}( V)$ of the matrix $H_B$. $H_B$ has a PHS if and only if $\text{Ad}_S(P)=-P$. Therefore, we define the set of relations over $\text{A}_{\text{loc}}$ 
$$\mathcal{R}_{\text{reduced}}=\{-\text{Ad}_S\}$$
This relation corresponds to an inner chiral symmetry in the language of \cite{kellendonk} and therefore (table \ref{final_classifications})
$$GL^{s.a.}_{h}(\mathbb{A}_{\text{loc}}, \mathcal{R})= GL^{s.a.}_{h}(\text{A}_{\text{loc}}, \mathcal{R}_{\text{reduced}})\leadsto DK(\text{A}_{\text{loc}}, \text{Ad}_S)\cong KU_1(\text{A}_{\text{loc}})$$

\subsection{Real cases without spin rotation symmetry}
\subsubsection{No symmetry}
In this case we just have the condition of being imaginary: therefore we only have the anti-symmetry with respect to the natural real structure $\text{Ad}_\gamma$.
$$\boxed{\mathcal{R}=\{-\text{Ad}_\gamma\}}$$
In \cite{kellendonk}, this symmetry corresponds to a PHS (anti-symmetry with respect to a real structure). Therefore (table \ref{final_classifications})

$$GL^{s.a.}_{h}(\mathbb{A}_{\text{loc}}, \mathcal{R})\leadsto KO_2(\mathbb{A}_{\text{loc}}^{\text{Ad}_\gamma})$$

\subsubsection{Time reversal symmetry}
$$\boxed{\mathcal{R}=\{-\text{Ad}_\gamma, \text{Ad}_{\Tilde{T}} \}}$$
These relations correspond to a PHS and TRS in the language of \cite{kellendonk}.

We define the product of the previous operators in order to obtain a chiral symmetry:
$$\alpha\coloneqq \text{Ad}_{\gamma}\circ\text{Ad}_{\Tilde{T}}=\text{Ad}_{\text{i}\gamma\Tilde{T}}=\text{Ad}_{\Gamma}$$
where $\Gamma\coloneqq \text{i}\gamma\Tilde{T}$ is the generator (self-adjoint and unitary) of the chiral symmetry. 
Therefore, we get the equivalent set of relations
$$\mathcal{R'}=\{\text{Ad}_{\Tilde{T}}, -\alpha \}$$
Moreover, we notice that 
$$\text{Ad}_{\Tilde{T}}(\Gamma)=-\Gamma$$
and so $\alpha$ is a imaginary inner grading. Therefore (paragraph \ref{chiralTRSSym}) 
$$GL^{s.a.}_{h}(\mathbb{A}_{\text{loc}}, \mathcal{R})=GL^{s.a.}_{h}(\mathbb{A}_{\text{loc}}, \mathcal{R}')\leadsto KO_{-1}(\mathbb{A}_{\text{loc}}^{\text{Ad}_{\Tilde{T}}})$$
However, in \cite{zirnbauer}, the natural (which we could also call \textit{reference}) real structure is $\text{Ad}_\gamma$.
We observe that $\mathfrak{f}=\text{Ad}_\gamma$ and $\mathfrak{r}=\text{Ad}_{\Tilde{T}}$ are \textit{inner related} real structures because

$$\mathfrak{f}\circ \mathfrak{r} =\text{Ad}_\Gamma$$
We compute the sign 
$$\eta^{(1)}_{\mathfrak{r},\mathfrak{f}}=\Gamma \mathfrak{r}(\Gamma)=-\text{i}\gamma\Tilde{T}\text{i}\gamma\Tilde{T}=-1$$
We are in the case of odd TRS with imaginary inner grading (table \ref{final_classifications}). Therefore,

$$GL^{s.a.}_{h}(\mathbb{A}_{\text{loc}}, \mathcal{R})\leadsto KO_{-1}(\mathbb{A}_{\text{loc}}^{\text{Ad}_{\Tilde{T}}})\cong KO_{3}(\mathbb{A}_{\text{loc}}^{\text{Ad}_\gamma})$$
To compare this case to the literature, we remark that this case corresponds to an odd TRS with even PHS.

\subsubsection{Time reversal symmetry + charge conservation}
$$\boxed{\mathcal{R}=\{-\text{Ad}_\gamma, \text{Ad}_{\Tilde{T}}, \text{Ad}_{Q} \}}$$
We can reduce and consider the algebra $\text{A}_{\text{loc}}$.
We define the set of relations over $\text{A}_{\text{loc}}$ 
$$\mathcal{R}_{\text{reduced}}=\{\text{Ad}_T\}$$
It corresponds to a TRS over $\text{A}_{\text{loc}}$. Then, table \ref{final_classifications} gives us
$$GL^{s.a.}_{h}(\mathbb{A}_{\text{loc}}, \mathcal{R})= GL^{s.a.}_{h}(\text{A}_{\text{loc}}, \mathcal{R}_{\text{reduced}})\leadsto KO_{0}(\text{A}_{\text{loc}}^{\text{Ad}_{T}}).$$

Using Lemma \ref{lemma1} we have
$$KO_{0}(\text{A}_{\text{loc}}^{\text{Ad}_{T}})\cong KO_{0}(\mathbb{A}_{\text{loc}}^{\text{Ad}_{\Tilde{T}}})$$
We want to compare to the natural real structure $\text{Ad}_{\gamma}$. As in the previous case we compute the relative sign for $\mathfrak{f}=\text{Ad}_\gamma$ and $\mathfrak{r}=\text{Ad}_{\Tilde{T}}$: $\eta^{(1)}_{\mathfrak{r},\mathfrak{f}}=-1$. Therefore we have an odd TRS and so (table \ref{final_classifications})
$$GL^{s.a.}_{h}(\mathbb{A}_{\text{loc}}, \mathcal{R})\leadsto KO_{4}(\mathbb{A}_{\text{loc}}^{\text{Ad}_{\gamma}})$$

\subsubsection{Time reversal symmetry + particle hole symmetry + charge conservation}
$$\boxed{\mathcal{R}=\{-\text{Ad}_\gamma, \text{Ad}_{\Tilde{T}}, \text{Ad}_{Q}, \text{Ad}_{L} \}}$$
We can reduce and we consider the algebra $\text{A}_{\text{loc}}$.
We define the set of relations over $\text{A}_{\text{loc}}$ 
$$\mathcal{R}_{\text{reduced}}=\{\text{Ad}_T,-\text{Ad}_S\}$$
As $\text{Ad}_T(S)=TST^{-1}=S$ and $S$ is unitary self-adjoint, $\mathcal{R}_{\text{reduced}}$ is equivalent to a TRS with a inner real chiral symmetry (table \ref{final_classifications}).
$$GL^{s.a.}_{h}(\mathbb{A}_{\text{loc}}, \mathcal{R})=GL^{s.a.}_{h}(\text{A}_{\text{loc}}, \mathcal{R}_{\text{reduced}})\leadsto  KO_{1}(\text{A}_{\text{loc}}^{\text{Ad}_{T}})\cong KO_{1}(\mathbb{A}_{\text{loc}}^{\text{Ad}_{\Tilde{T}}})$$

We compare as before to the reference real structure $\text{Ad}_{\gamma}$ and we get (table \ref{final_classifications}, odd TRS inner real chiral symmetry)

$$GL^{s.a.}_{h}(\mathbb{A}_{\text{loc}}, \mathcal{R})\leadsto KO_{5}(\mathbb{A}_{\text{loc}}^{\text{Ad}_{\gamma}})$$
To compare this case to the literature, we should express this symmetry in terms of TRS and PHS. If we compose the TRS real structure and the chiral symmetry we can find the corresponding PHS real structure $\mathfrak{p}$
$$\mathfrak{p}=\text{Ad}_T\circ\text{Ad}_S=\text{Ad}_{TS}$$
As usual $\mathfrak{f}=\text{Ad}_\gamma$ is the reference real structure. A generator for $\mathfrak{f}\circ \mathfrak{p}$ is $u=\gamma T S$. Therefore the sign $\eta^{(1)}_{\mathfrak{p},\mathfrak{f}}$ is 
$$\eta^{(1)}_{\mathfrak{p},\mathfrak{f}}=u\mathfrak{p}(u)=-1.$$
Therefore we are in the case of an odd TRS with odd PHS.

\subsection{Real cases with spin rotations symmetry}

The SRS is a remarkable symmetry because we will see that the commutation with the spin rotation operators allows us to factorize the algebra using quaternions. For each one of the following cases, we will show the existence of a real isomorphism of the type
$$(\mathcal{L}( W), \text{Ad}_\gamma)\cong (\mathcal{L}( W^+)\otimes \quat,\mathfrak{r}\otimes \text{id})$$
for a particular (smaller) Real algebra $(\mathcal{L}( W^+), \mathfrak{r})$.

On the other hand, we know that we can describe quaternions $\quat$ as 
$$\quat\cong\text{span}_\R \{1,\text{i}\sigma_1,\text{i}\sigma_2,\text{i}\sigma_3 \}$$
The role of quaternions when dealing with real $K$-theory is made explicit with the help of the following, well known, result.
\begin{theorem}\label{lemma2}
    Let $B$ a real $C^*$-algebra. Then for all $i\in \{0, \dots,7\}$
    $$KO_i(B\otimes \quat)\cong KO_{i+4}(B)$$
    ($i+4$ is taken modulo $8$).
\end{theorem}
\begin{proof}
    Equation (\ref{KOi}) gives us
    $$KO_{s-r+1}(B)\cong DK(B\otimes Cl_{r,s}, \text{id}\otimes \text{st})$$
    Then, using the properties of graded Clifford algebras
    \begin{align*}
        KO_{i} (B\otimes \quat) &= DK(B\otimes \quat \otimes Cl_{2,i+1}, \text{id}\otimes \text{id}\otimes\text{st})\quad & \text{by definition}\\
        &=DK(B\otimes \quat \otimes (Cl_{1,1}\hat{\otimes} Cl_{1,i}), \text{id}\otimes \text{id}\otimes\text{st}) \quad & \text{eq. (\ref{Cl1})}\\
        &=DK(B\otimes Cl_{0,4}\hat{\otimes} Cl_{1,i}, \text{id}\otimes \text{st}\otimes\text{st}) \quad & \text{eq. (\ref{Cl2})}\\
        &= DK(B\otimes Cl_{1,i+4}, \text{id}\otimes \text{st}) \quad & \text{(eq. \ref{Cl1})}\\
        &= KO_{i+4} (B)\quad & \text{by definition}
    \end{align*}
\end{proof}
The advantage of the previous factorisation is that, with our choice for $(\mathcal{L}( W^+), \mathfrak{r})$, the relations of commutation with $\Tilde{j}_\mu$ translate in a more readable way and they are easier to handle. 
In the following lemma (which is in the spirit very similar to the lemma 4.1.26 of \cite{CMax}) we sum up the previous considerations and we describe with more details the construction of the real algebra and of the isomorphism.
\begin{theorem}\label{lemma3}
    Given the generators of spin rotations $\Tilde{j}_1$, $\Tilde{j}_2$ and $\Tilde{j}_3$, we define $ W^\pm\coloneqq \text{ker}(\Tilde{j}_3\mp\text{i})$, the anti-linear operator $\gamma^+\coloneqq-\gamma {\Tilde{j}_1}\vert_{ W^+}$ and $\mathbb{A}_{\text{loc}}^+\coloneqq \mathcal{L}( W^+)$.
    Then we have a real isomorphism 
    $\mathbb{A}_{\text{loc}}^{\text{Ad}_\gamma} \cong (\mathbb{A}_{\text{loc}}^+)^{\text{Ad}_{\gamma^+}}\otimes \quat$.
    Moreover, if we define the relations $\mathcal{R}=\{-\text{Ad}_\gamma,\text{Ad}_{\Tilde{j}_1},\text{Ad}_{\Tilde{j}_2} \}$ over $\mathbb{A}_{\text{loc}}$ and $\mathcal{R}'=\{-\text{Ad}_{\gamma^+} \}$ over $\mathbb{A}_{\text{loc}}^+$, we have a homeomorphism of topological spaces between $GL^{s.a.}(\mathbb{A}_{\text{loc}}, \mathcal{R})$ and $ GL^{s.a.}(\mathbb{A}_{\text{loc}}^+, \mathcal{R}')$ (where $GL^{s.a.}(A, \mathcal{R})$ is the set of self-adjoint invertible elements of $A$ satisfying the relations of $\mathcal{R}$).
\end{theorem}
\begin{proof}
    As $\Tilde{j}_3$ is skew-hermitian and unitary, we have that $\pm$i are the only eigenvalues for $\Tilde{j}_3$. Therefore we have the decomposition $ W= W^+\oplus  W^-$. With respect to this decomposition, we can write the $\Tilde{j}_\mu$'s and $\gamma$ as
$$\Tilde{j}_1=\begin{pmatrix}
0 & \Tilde{j}_1\vert_{ W^-}\\
\Tilde{j}_1\vert_{ W^+} & 0
\end{pmatrix}\quad
\Tilde{j}_2=\begin{pmatrix}
0 & \text{i}\Tilde{j}_1\vert_{ W^-}\\
-\text{i}\Tilde{j}_1\vert_{ W^+} & 0
\end{pmatrix}
\quad
\Tilde{j}_3=\begin{pmatrix}
\text{i} & 0\\
0 & -\text{i}
\end{pmatrix}
\quad
\gamma=\begin{pmatrix}
0 & \gamma\vert_{ W^-}\\
\gamma\vert_{ W^+} & 0
\end{pmatrix}$$

We define the map $\phi$
$$\phi\coloneqq \begin{pmatrix}
1 & 0\\
0 & -\Tilde{j}_1\vert_{ W^+}
\end{pmatrix}: W^+\oplus W^+\rightarrow  W^+\oplus W^-$$
which is invertible and the inverse is given by
$$\phi^{-1}=\begin{pmatrix}
1 & 0\\
0 & \Tilde{j}_1\vert_{ W^-}
\end{pmatrix}$$
We would like to find a real structure $\mathfrak{r}=\gamma^+\otimes \mathfrak{t}$ over $ W^+\otimes \C^2$ such that $\phi^{-1}$ is a real isomorphism between $( W, \gamma)$ and $( W^+\otimes \C^2, \gamma^+\otimes \mathfrak{t})$, where 
$$\mathfrak{t}=\begin{pmatrix}
0 & \mathfrak{c} \\
-\mathfrak{c} & 0
\end{pmatrix}=\text{i}\sigma_y\circ\mathfrak{c}$$
$\phi^{-1}$ is real if $\phi^{-1}\circ \gamma=\mathfrak{r}\circ\phi^{-1}$. Therefore (we use the fact that $\Tilde{j}_1$ is real, namely it commutes with $\gamma$)
$$\mathfrak{r}= \phi^{-1}\gamma\phi= \begin{pmatrix}
	0 & -\Tilde{j}_1\gamma\vert_{ W^+}\\
	\Tilde{j}_1\gamma\vert_{ W^+} & 0
	\end{pmatrix}= \gamma^+\otimes \mathfrak{t}$$
where $\gamma^+\coloneqq -\Tilde{j}_1\gamma$.
The previous isomorphism induces the real isomorphism 
\begin{align*}
\chi:  (\mathcal{L}( W), \text{Ad}_\gamma)&\rightarrow (\mathcal{L}( W^+)\otimes M_2(\mathbb{C}), \text{Ad}_{\gamma^+\otimes \mathfrak{t}})\\  x &\mapsto \text{Ad}_{\phi^{-1}}(x)
\end{align*}
As $M_2(\mathbb{C})^{\text{Ad}_{\mathfrak{t}}}\cong \quat$, we have that $\chi$ is an isomorphism between $\mathbb{A}_{\text{loc}}^{\text{Ad}_\gamma} $ and $ (\mathbb{A}_{\text{loc}}^+)^{\text{Ad}_{\gamma^+}}\otimes \quat$.

We show now the bijection between $GL^{s.a.}(\mathbb{A}_{\text{loc}}, \mathcal{R})$ and $GL^{s.a.}(\mathbb{A}_{\text{loc}}^+, \mathcal{R}')$. First of all, we verify that
$$\chi(\Tilde{j}_1)=1\otimes \text{i}\sigma_y\quad \text{and}\quad \chi(\Tilde{j}_2)=1\otimes \text{i}\sigma_x$$ If we take $H\in GL^{s.a.}(\mathbb{A}_{\text{loc}}, \mathcal{R})$, we have ($\chi $ preserves the real structure, anti-commuting relations and the unity)
\begin{itemize}
	\item $\chi(H)$ is imaginary. Therefore
	$$\chi(H)=z_0\otimes \text{i}\sigma_x+z_1\otimes \text{i}\sigma_y+z_2\otimes \text{i}\sigma_z+z_3\otimes 1$$
	with $z_i$ imaginary (because imaginary elements of $\mathcal{L}( W^+)\otimes M_2(\mathbb{C})$ are $\R$-linear combinations of elements of the form $a\otimes b $ with $a$ imaginary and $b$ real and we have a basis for the real elements of $M_2(\mathbb{C})$).
	\item $[\chi(H),1\otimes \text{i}\sigma_y]=0$. Then
	$$\chi(H)=z_1\otimes \text{i}\sigma_y+z_3\otimes 1$$
	\item $[\chi(H),1\otimes \text{i}\sigma_x]=0$. Then
	$$\chi(H)=z_3\otimes 1$$
	\item $\chi(H)$ is self-adjoint. Therefore $z_3$ is self-adjoint.
	\item $\chi(H)$ is unitary. Therefore $z_3$ is a unitary.    
\end{itemize}
\textbf{Remark:} $H$ commuting with $\Tilde{j}_3$ does not say anything more: $\Tilde{j}_3=\Tilde{j}_2 \Tilde{j}_1$!

To sum up, $H\in GL^{s.a.}(\mathbb{A}_{\text{loc}}, \mathcal{R})$ if and only if $\chi(H)=z\otimes 1$ with $z$ unitary imaginary self-adjoint element of $\mathcal{L}( W^+)$ and we get the isomorphism we wanted.  
\end{proof}

Below we will note, with a slight abuse of notation, $\mathfrak{f}$ the real structure $\text{Ad}_{\gamma}$ over $\mathbb{A}_{\text{loc}}$ or $\text{Ad}_{\gamma^+}\otimes \text{Ad}_{\mathfrak{t}}$ over $\mathbb{A}_{\text{loc}}^+\otimes M_2(\C)$.

\subsubsection{Only spin rotations symmetry}
$$\boxed{\mathcal{R}=\{-\text{Ad}_\gamma,\text{Ad}_{\Tilde{j}_1},\text{Ad}_{\Tilde{j}_2} \}}$$

If we set $\mathcal{R}'=\{-\text{Ad}_{\gamma^+} \}$, lemma \ref{lemma3} gives us
$$GL^{s.a.}_{h}(\mathbb{A}_{\text{loc}}, \mathcal{R})\cong GL^{s.a.}_{h}(\mathbb{A}_{\text{loc}}^+, \mathcal{R}').$$
We are in the case of a PHS and table \ref{final_classifications} gives us
$$GL^{s.a.}_{h}(\mathbb{A}_{\text{loc}}^+, \mathcal{R}')\leadsto KO_{2}((\mathbb{A}_{\text{loc}}^+)^{\text{Ad}_{\gamma^+}})$$
As $K$-theory is stable under Morita equivalence, we know that 
$$KO_{2}((\mathbb{A}_{\text{loc}}^+)^{\text{Ad}_{\gamma^+}})\cong KO_{2}((\mathbb{A}_{\text{loc}}^+\otimes M_2(\C))^{\text{Ad}_{\gamma^+}\otimes \mathfrak{c}})$$
which is the $K$-group corresponding to a PHS over the algebra $\mathbb{A}_{\text{loc}}^+\otimes M_2(\C)$ with real structure $\mathfrak{p}=\text{Ad}_{\gamma^+}\otimes \mathfrak{c}$. 
If we take the reference real structure $\mathfrak{f}$ defined as above, we have 
$$\mathfrak{f}\circ \mathfrak{r}=\text{Ad}_{\text{id}\otimes \text{i}\sigma_y}$$
and we get the relative sign of $\mathfrak{p}$ with respect to $\mathfrak{f}$ $$\eta^{(1)}_{\mathfrak{p},\mathfrak{f}}=(\text{id}\otimes \text{i}\sigma_y )\mathfrak{p}(\text{id}\otimes \text{i}\sigma_y)=-1.$$
Therefore we have an odd PHS and comparing to table \ref{final_classifications} we get  
$$GL^{s.a.}_{h}(\mathbb{A}_{\text{loc}}, \mathcal{R})\leadsto   KO_6(\mathbb{A}_{\text{loc}}^{\text{Ad}_\gamma}).$$

\subsubsection{SRS + time reversal symmetry}
$$\boxed{\mathcal{R}=\{-\text{Ad}_\gamma,\text{Ad}_{\Tilde{j}_1},\text{Ad}_{\Tilde{j}_2},\text{Ad}_{\Tilde{T}} \}}$$
We define the product of $\text{Ad}_\gamma$ and $\text{Ad}_{\Tilde{T}}$ in order to obtain a chiral symmetry:
$$\alpha\coloneqq \text{Ad}_{\gamma}\circ\text{Ad}_{\Tilde{T}}=\text{Ad}_{\text{i}\gamma\Tilde{T}}=\text{Ad}_{\Gamma}$$
where $\Gamma\coloneqq \text{i}\gamma\Tilde{T}$ is the generator (self-adjoint and unitary) of the chiral symmetry. 

We consider the decomposition $ W= W^+\oplus  W^-$ as in lemma \ref{lemma3}. As $\Tilde{j}_3$ commutes with $\Tilde{T}$, we have that if $\psi\in  W^\pm$ then
$$\Tilde{j}_3 \Tilde{T}\psi=\Tilde{T}\Tilde{j}_3 \psi=\mp\text{i}\Tilde{T}\psi \iff \Tilde{T}\psi\in  W^\mp$$
Therefore, with respect to $ W= W^+\oplus  W^-$
$$\Tilde{T}=\begin{pmatrix}
0 & \Tilde{T}\vert_{ W^-}\\
\Tilde{T}\vert_{ W^+} & 0
\end{pmatrix}$$
We consider the isomorphism $\chi$ of lemma \ref{lemma3} and we compute the image $\xi\in \mathcal{L}( W^+)\otimes M_2(\mathbb{C})$ of the generator of the chiral symmetry $\Gamma$ through $\chi$.
\begin{align*}
    \xi&=\chi(\Gamma)=\phi^{-1}\text{i}\gamma \Tilde{T}\phi\\&=
    \begin{pmatrix}
        1 & 0\\
        0 & \Tilde{j}_1\vert_{ W^-}         
    \end{pmatrix}
    \begin{pmatrix}
        0 & \text{i}\gamma\vert_{ W^-}\\
        \text{i}\gamma\vert_{ W^+} & 0        
    \end{pmatrix}
    \begin{pmatrix}
        0 & \Tilde{T}\vert_{ W^-}\\
        \Tilde{T}\vert_{ W^+} & 0        
    \end{pmatrix}
    \begin{pmatrix}
        1 & 0\\
        0 & -\Tilde{j}_1\vert_{ W^+}         
    \end{pmatrix}\\&=
    \text{i}\gamma \Tilde{T}\otimes \text{id}
\end{align*}
Again, by stability of $K$-theory by Morita equivalence, we can consider the natural extension of the relations $\{-\text{Ad}_{\gamma^+} , -\text{Ad}_{\text{i}\gamma \Tilde{T}}\}$ to $\mathbb{A}_{\text{loc}}^+\otimes M_2(\C)$: we set $\mathcal{R}'=\{-\mathfrak{p} , -\alpha\}$, where $\mathfrak{p}=\text{Ad}_{\gamma^+}\otimes \mathfrak{c} $ and $\alpha=\text{Ad}_{\text{i}\gamma \Tilde{T}\otimes \text{id}}$.
$\mathcal{R}'$ corresponds to a PHS and a chiral symmetry on $\mathbb{A}_{\text{loc}}^+\otimes M_2(\C)$. We take the product of $\mathfrak{p}$ and $\alpha$ to find the corresponding TRS.
$$\mathfrak{r}=\mathfrak{p}\circ \alpha=\text{Ad}_{\text{i}\Tilde{j}_1 \Tilde{T}}\otimes\mathfrak{c}$$
We define $\mathcal{R}''=\{\mathfrak{r}, -\alpha\}$.
We compute
$$\mathfrak{r}(\text{i}\gamma \Tilde{T}\otimes\text{id})=-\text{i}\gamma \Tilde{T}\otimes\text{id}$$
and therefore we have a imaginary inner grading. 
Taking as usual $\mathfrak{f}=\text{Ad}_{\gamma^+}\otimes \text{Ad}_{\mathfrak{t}}$ as reference real structure, we have that $u=\text{i}\gamma\Tilde{T}\otimes \text{i}\sigma_y$ is a generator of $\mathfrak{f}\circ \mathfrak{r}$ and we get the relative sign
$$\eta^{(1)}_{\mathfrak{r},\mathfrak{f}}=u \mathfrak{r}(u)=+1$$
We compare with table \ref{final_classifications}, for the case imaginary inner grading with even TRS. Therefore
$$GL^{s.a.}_{h}(\mathbb{A}_{\text{loc}}^+\otimes M_2(\C), \mathcal{R}'')\leadsto KO_{7}(\mathbb{A}_{\text{loc}}^{\text{Ad}_\gamma})$$
To sum up
$$GL^{s.a.}_{h}(\mathbb{A}_{\text{loc}}, \mathcal{R})\leadsto KO_7(\mathbb{A}_{\text{loc}}^{\text{Ad}_\gamma})$$
To compare this case to the literature, we can verify that this case corresponds to an even TRS ($\mathfrak{r}$) with odd PHS ($\mathfrak{p}$) over $\mathbb{A}_{\text{loc}}^+\otimes M_2(\C)$ with respect to $\mathfrak{f}$.

\subsubsection{SRS + time reversal symmetry + charge conservation}
$$\boxed{\mathcal{R}=\{-\text{Ad}_\gamma,\text{Ad}_{\Tilde{j}_1},\text{Ad}_{\Tilde{j}_2},\text{Ad}_{\Tilde{T}},\text{Ad}_Q  \}}$$
We can reduce and we consider the algebra $\text{A}_{\text{loc}}$.
We define the set of relations over $\text{A}_{\text{loc}}$ 
$$\mathcal{R}_{\text{reduced}}=\{\text{Ad}_T,\text{Ad}_{j_1},\text{Ad}_{j_2}\}.$$

We state here a result very similar to lemma \ref{lemma3}.
\begin{theorem}\label{lemma4}
    Given the generators of spin rotations $j_1$, $j_2$ and $j_3$ and the time-reversal operator $T$ over $ V$, we define $ V^\pm\coloneqq \text{ker}(j_3\mp\text{i})$, the anti-linear operator $T^+\coloneqq-j_1 T\vert_{ V^+}$ and $\text{A}_{\text{loc}}^+\coloneqq \mathcal{L}( V^+)$.
    Then we have a real isomorphism 
    $\text{A}_{\text{loc}}^{\text{Ad}_T} \cong (\text{A}_{\text{loc}}^+)^{\text{Ad}_{T^+}}\otimes \quat$.
    Moreover, if we define the relations $\mathcal{R}_{\text{reduced}}=\{\text{Ad}_T,\text{Ad}_{j_1},\text{Ad}_{j_2} \}$ over $\text{A}_{\text{loc}}$ and $\mathcal{R}'_{\text{reduced}}=\{\text{Ad}_{T^+} \}$ over $\text{A}_{\text{loc}}^+$, we have a homeomorphism of topological spaces 
    $$GL^{s.a.}_{h}(\text{A}_{\text{loc}}, \mathcal{R}_{\text{reduced}})= GL^{s.a.}_{h}(\text{A}_{\text{loc}}^+, \mathcal{R}'_{\text{reduced}})$$
\end{theorem}
\begin{proof}
    The proof is the same as for lemma \ref{lemma3} where we consider as isomorphism $\phi$ 
    $$\phi\coloneqq \begin{pmatrix}
        1 & 0\\
        0 & -j_1\vert_{ V^+}
    \end{pmatrix}: V^+\oplus V^+\rightarrow  V^+\oplus V^-$$
    We find the expression for $T^+$ computing the image $\Gamma$ via $\text{Ad}_{\phi^{-1}}$ of $T$.
    \begin{align*}
    \Gamma&=\phi^{-1}T\phi\\&=
    \begin{pmatrix}
        1 & 0\\
        0 & j_1\vert_{ V^-}         
    \end{pmatrix}
    \begin{pmatrix}
        0 & T\vert_{ V^-}\\
        T\vert_{ V^+} & 0        
    \end{pmatrix}
    \begin{pmatrix}
        1 & 0\\
        0 & -j_1\vert_{ V^+}         
    \end{pmatrix}=
    \begin{pmatrix}
        0 & -j_1 T\vert_{ V^-}\\
        j_1T\vert_{ V^+} & 0        
    \end{pmatrix}\\&=-j_1 T\otimes \begin{pmatrix}
        0 & \mathfrak{c}\\
        -\mathfrak{c} & 0        
    \end{pmatrix}=
    T^+\otimes \mathfrak{t}
\end{align*}
The previous isomorphism induces the real isomorphism 
\begin{align*}
    \chi:  (\mathcal{L}( V), \text{Ad}_T)&\rightarrow (\mathcal{L}( V^+)\otimes M_2(\mathbb{C}), \text{Ad}_{T^+\otimes \mathfrak{t}})\\  x &\mapsto \text{Ad}_{\phi^{-1}}(x)
\end{align*}
\end{proof}

Lemma \ref{lemma4} tells us that we can just consider the relation $\mathcal{R}'_{\text{reduced}}=\{\text{Ad}_{T^+}\}$ over $\text{A}_{\text{loc}}^+$ and we have
$$GL^{s.a.}_{h}(\text{A}_{\text{loc}}, \mathcal{R}_{\text{reduced}})= GL^{s.a.}_{h}(\text{A}_{\text{loc}}^+, \mathcal{R}'_{\text{reduced}}).$$
The relation of $\mathcal{R}'_{\text{reduced}}$ corresponds to a TRS and we would like to compare it to the reference real structure $\mathfrak{f}=\text{Ad}_{\gamma}$ over $\mathbb{A}_{\text{loc}}$. To this end, by stability of $K$-theory by Morita equivalence, we can consider the natural extension $\mathfrak{r}$ of the relation of $\mathcal{R}'_{\text{reduced}}$ to $\mathbb{A}_{\text{loc}}$. The isomorphisms $\psi$ of presented in (\ref{M2(A)}) and $\chi$ defined in the proof of lemma \ref{lemma4} allow us to set $\mathfrak{r}=\text{Ad}_{\Gamma}$ where 
\begin{align*}
    \Gamma&=\psi((\chi^{-1}(T^+\otimes \mathfrak{c}))\otimes \mathfrak{c}) \\&=\mathfrak{M}((\phi(T^+\otimes \mathfrak{c})\phi^{-1})\otimes \mathfrak{c})\mathfrak{M}^{-1}\\&=\begin{pmatrix}
    1 & 0\\
    0 & \mathfrak{m}
\end{pmatrix}
    \begin{pmatrix}
        1 & 0\\
        0 & -j_1\vert_{ V^+}         
    \end{pmatrix}
    \begin{pmatrix}
        T^+\vert_{ V^+} & 0\\
        0 & T^+\vert_{ V^+}        
    \end{pmatrix}
    \begin{pmatrix}
        1 & 0\\
        0 & j_1\vert_{ V^-}         
    \end{pmatrix}\begin{pmatrix}
    1 & 0\\
    0 & \mathfrak{m}^{-1}
\end{pmatrix}\\&=
    -\Tilde{j}_1 \Tilde{T}.
\end{align*}
We see that
$$\mathfrak{f}\circ \mathfrak{r}=\text{Ad}_{\gamma}\circ \text{Ad}_{\Gamma}=\text{Ad}_{-\gamma\Tilde{j}_1 \Tilde{T}}$$
and therefore $u=\gamma\Tilde{j}_1 \Tilde{T}$ is a generator for $\mathfrak{f}\circ \mathfrak{r}$. We have the relative sign
$$\eta^{(1)}_{\mathfrak{r},\mathfrak{f}}=u\mathfrak{r}(u) =+1.$$
We are in the case of an even TRS and therefore table \ref{final_classifications} gives us
$$GL^{s.a.}_{h}(\mathbb{A}_{\text{loc}}, \mathcal{R})\leadsto KO_{0}(\mathbb{A}_{\text{loc}}^{\text{Ad}_\gamma}).$$

\subsubsection{SRS + time reversal symmetry + charge conservation + particle hole symmetry}
$$\boxed{\mathcal{R}=\{-\text{Ad}_\gamma,\text{Ad}_{\Tilde{j}_1},\text{Ad}_{\Tilde{j}_2},\text{Ad}_{\Tilde{T}},\text{Ad}_Q,-\text{Ad}_L  \}}$$
We can reduce and we consider the algebra $\text{A}_{\text{loc}}$.
We define the set of relations over $\text{A}_{\text{loc}}$ 
$$\mathcal{R}_{\text{reduced}}=\{\text{Ad}_T,-\text{Ad}_S,\text{Ad}_{j_1},\text{Ad}_{j_2}\}$$
We apply lemma \ref{lemma4}. As $S$ is a $\C$-linear operator commuting with $j_3$, we can express it with respect to $ V\cong  V^+\oplus V^-$ as
$$S=\begin{pmatrix}
    S\vert_{ V^+} & 0\\
    0 & S\vert_{ V^-}
\end{pmatrix}$$
We compute the operator corresponding to the chiral symmetry over $(\text{A}_{\text{loc}}^+\otimes M_2(\mathbb{C}))$.
$$L\coloneqq\phi^{-1}S\phi=
    \begin{pmatrix}
    S\vert_{ V^+} & 0\\
    0 & S\vert_{ V^+}
\end{pmatrix}
    =S\vert_{ V^+}\otimes\text{id}$$
We define the set of relations over $\text{A}_{\text{loc}}^+$ 
$$\mathcal{R}'_{\text{reduced}}=\{\text{Ad}_{T^+},-\text{Ad}_S\}$$
Therefore we have a TRS with inner real chiral symmetry because
$$\text{Ad}_{T^+}(S)=j_1TSj_1T=S$$
As in the previous case, $\mathfrak{r}=\text{Ad}_{\Tilde{j}_1 \Tilde{T}}$ is the natural extension of $\text{Ad}_{T^+}$ to $\mathbb{A}_{\text{loc}}$ and the relative sign with respect to the reference real structure $\mathfrak{f}=\text{Ad}_{\gamma}$ is
$$\eta^{(1)}_{\mathfrak{r},\mathfrak{f}}=+1.$$
We are in the case of an even TRS with inner real chiral symmetry and therefore table \ref{final_classifications} gives us
$$GL^{s.a.}_{h}(\mathbb{A}_{\text{loc}}, \mathcal{R})\leadsto KO_{1}(\mathbb{A}_{\text{loc}}^{\text{Ad}_\gamma}).$$

We can verify that if $\mathfrak{p}$ is the natural extension of the real structure obtained composing $\text{Ad}_T$ and $\text{Ad}_S $ we obtain an even PHS. Therefore we are in the case of an even TRS with odd PHS.
\subsection{To sum up}
Given the algebra of bulk observables and charge-conserving observables $\mathbb{A}$ and A (containing the information of the internal and spacial degrees of freedom), we can summarise in a table what we described in the previous paragraph.
$$\begin{array}{c c | c c |c| c c c}
 & &\text{\cite{zirnbauer}} & & &\text{\cite{kellendonk}}\\
\hline 
 s & \text{Class} & & \text {symmetries} & K\textbf {-group}&\text{symmetries}&\text {sign TRS} &\text{sign PHS}\\
\hline 
0 & \text{D} & \text{Real cases} &  \text {none} & KO_2(\mathbb{A}^{\text{Ad}_{\gamma}})&  \text{PHS} & & +1\\
1 & \text{DIII} & \text{without SRS} &  \text {TRS} & KO_{3}(\mathbb{A}^{\text{Ad}_{\gamma}})&  \text{TRS+PHS}& -1 & +1\\
2 & \text{AII} &  &  \text {TRS+Q} & KO_{4}(\mathbb{A}^{\text{Ad}_{\gamma}})& \text{TRS}& -1\\
3 & \text{CII} &  &  \text {TRS+Q+PHS} & KO_{5}(\mathbb{A}^{\text{Ad}_{\gamma}})&\text{TRS+PHS} & -1 & -1\\
\hline
4 & \text{C} & \text{Real cases} & \text {none} & KO_{6}(\mathbb{A}^{\text{Ad}_{\gamma}})&\text{PHS}  & & -1\\
5 & \text{CI} & \textbf{with SRS} &  \text {TRS} & KO_{7}(\mathbb{A}^{\text{Ad}_{\gamma}})&\text{TRS+PHS}  & +1 & -1\\
6 & \text{AI} & &  \text {TRS+Q} & KO_{0}(\mathbb{A}^{\text{Ad}_{\gamma}})&  \text{TRS} & +1 &\\
7 & \text{BDI} & &  \text {TRS+Q+PHS} & KO_{1}(\mathbb{A}^{\text{Ad}_{\gamma}})&\text{TRS+PHS} & +1 & +1\\
\hline
0 & \text{A} & \text{Complex cases}  &  \text { Q } &  KU_0(\text{A}) &  \text{none} \\
1 & \text{AIII} & & \text { Q+PHS } & KU_1(\text{A})&  \text{chiral}\\
\end{array}
$$

The first column corresponds to the number of symmetries, while the "Class" column corresponds to the standard name in the Kitaev table. What we obtained is a translation of the symmetries in the concrete language of \cite{zirnbauer} (left part of the table) into the abstract language of \cite{kellendonk} (right part). The central column shows which is the $K$-group obtained following one of the two approaches (we are glad to see that both approaches lead to the same $K$-group!). For example, we see that a Time Reversal Symmetry without Spin Rotations Symmetry in \cite{zirnbauer} correspond do an odd Time Reversal Symmetry with an even Particle Hole Symmetry in \cite{kellendonk} and both lead to the $K$-group $KO_{3}(\mathbb{A}^{\text{Ad}_{\gamma}})$.

We conclude this work summarising the main remarks about the comparison of the two approaches which are emphasised by the table.
\begin{itemize}
    \item \textbf{The reference real structure.} In order to compare different $K$-groups it is useful to choose a reference real structure. Indeed, \cite{kellendonk} shows how we can relate the $K$-groups coming from two Real $C^*$-algebras when the underlying complex algebra is the same and we know the relative sign between the two real structures. On the other hand, we saw that the real structure $\gamma$ appeared naturally in the construction of the Nambu space. Therefore, it is reasonable to consider the induced real structure $\text{Ad}_\gamma$ as a reference real structure. In this way, the eight real $K$-groups have the same underlying Real $C^*$-algebra and it makes sense to compare them.
    \item \textbf{Charge conservation.} The symmetry of charge conservation (Q) annihilates the real symmetry $\text{Ad}_\gamma$. This means that, imposing charge conservation, we obtain operators which are a double copy of operators over the complex algebra A and which satisfy trivially the natural real structure. Therefore, in the case where no other real symmetries are present, $Q$ allows us to reduce the Real $C^*$-algebra $\mathbb{A}^{\text{Ad}_{\gamma}}$ and just consider the complex algebra A. This leads to the last two lines of the table, where we obtain the two complex $K$-groups. 
    \item \textbf{Spin Rotations Symmetry.} The Spin Rotations Symmetry (SRS) allows us to factor the Real $C^*$-algebra $\mathbb{A}^{\text{Ad}_{\gamma}}$ using quaternions: 
    $$\mathbb{A}^{\text{Ad}_\gamma} \cong (\mathbb{A}^+)^{\text{Ad}_{\gamma^+}}\otimes \quat.$$ 
    Therefore, we can consider a fewer number of relations over the smaller algebra $\mathbb{A}^+$. 
    This leads to a change of the sign of the real structures with respect to the reference real structure $\text{Ad}_\gamma$ in \cite{kellendonk}. In terms of $K$-groups this fact shows up in the shifting of 4 units of the $K$-group indices.  
    \item \textbf{Odd Time Reversal Symmetry.} We observe that in the first rows of the table the sign of the Time Reversal Symmetry is always $-1$. This is due to the fact that we are dealing with fermions and therefore the time reversal operator squares to $-1$. We need Spin Rotation Symmetry to flip the relative sign.
\end{itemize}

\bibliographystyle{alpha}
\bibliography{bibliography} 
\end{document}